\documentclass[12pt]{amsart}

\UseRawInputEncoding

\newtheorem{theorem}{Theorem}[section]
\newtheorem{lemma}[theorem]{Lemma}
\newtheorem{corollary}[theorem]{Corollary}

\newtheorem{conjecture}[theorem]{Conjecture}

\theoremstyle{plain}
\newtheorem{definition}[theorem]{Definition}

\newtheorem{question}[theorem]{Question}

\newtheorem{setting}[theorem]{Setting}

\theoremstyle{definition}

\theoremstyle{remark}

\numberwithin{equation}{section}

\usepackage{fancyhdr}
\usepackage{amscd}
\usepackage{amsmath}
\usepackage{amsthm}
\usepackage{amssymb}

\begin{document}

\title[K-theory and obstructions]{Bloch-Ogus theorem, cyclic homology and deformation of Chow groups}

\author{Sen Yang }

\address{School of Mathematics \\  Southeast University \\
Nanjing, China\\
}

\address{Shing-Tung Yau Center of Southeast University \\ 
Southeast University \\
Nanjing, China\\
}

\email{101012424@seu.edu.cn}

\date{}

\maketitle

\begin{abstract}
Using Bloch-Ogus theorem and Chern character from K-theory to cyclic homology, we answer a question of Green and Griffiths on extending Bloch formula. Moreover, we construct a map from local Hilbert functor to local cohomology. With suitable assumptions, we use this map to answer a question of Bloch on constructing a natural transformation from local Hilbert functor to cohomological Chow groups.
\end{abstract}


\section{Introduction}
This paper is devoted to studying infinitesimal deformation of Chow group $CH^{p}(X)$ of codimension $p$ algebraic cycles modulo rational equivalence, where $X$ is a smooth projective variety over a field $k$ of characteristic zero. After many years' intensive study, the structure of $CH^{p}(X)$ for general $p$ still remains largely open. To understand  Chow groups infinitesimally, Bloch pioneered to study formal completions of $CH^{p}(X)$. One fundamental tool in this approach is Bloch formula (cf. Bloch \cite{Bl2-Annals}, Quillen \cite{Quillen} and Soul\'e \cite{Soule}),
{\footnotesize
\begin{equation}
CH^{p}(X)_{\mathbb{Q}}=H^{p}(X, K^{M}_{p}(O_{X}))_{\mathbb{Q}},
\end{equation}
}where $K^{M}_{p}(O_{X})$ is the Milnor K-theory sheaf associated to the presheaf $U \to K^{M}_{p}(O_{X}(U))$ with $U \subset X$ open affine. Kerz generalized the isomorphism (1.1) in \cite{Kerz}.

Bloch formula motivates two functors on the category $Art_{k}$ (see Notation (2) on page 4 below)
{\footnotesize
\begin{align}
\widetilde{CH}^{p}: & \ A \to H^{p}(X, K^{M}_{p}(O_{X_A}))_{\mathbb{Q}}, \\
\widehat{CH}^{p}:  & \ A \to \mathrm{kernel \ of } \ \{ H^{p}(X, K^{M}_{p}(O_{X_A}))_{\mathbb{Q}}  \xrightarrow{aug} H^{p}(X, K^{M}_{p}(O_{X}))_{\mathbb{Q}} \},
\end{align}
}where $A \in Art_{k}$, $X_A=X\times_{\mathrm{Spec}(k)}\mathrm{Spec}(A)$ and $aug$ is the map induced by augmentation $A \to k$. The group $\widetilde{CH}^{p}(A)$, which is called cohomological Chow group, can be considered as deformation of Chow group $CH^{p}(X)$, and $\widehat{CH}^{p}(A)$ is called formal completion of $CH^{p}(X)$, see Bloch \cite{Bl3} and Stienstra \cite{Stien1}. These functors $\widetilde{CH}^{p}$ and $\widehat{CH}^{p}$ are closely related with some major conjectures, including variational (infinitesimal) Hodge conjecture. We refer to Bloch-Esnault-Kerz \cite{BEK1, BEK2}, Green-Griffiths \cite{GGChow}, Morrow \cite{Morrow} and Patel-Ravindra \cite{PatelRavi1} for recent progress on these conjectures.

Bloch \cite{Bl3} studied these two functors in the case that $k$ is a number field and asked the following important conjecture.
\begin{conjecture} [\cite{Bl3}] \label{c:Bloch}
Let $X$ be a smooth complex projective surface with trivial geometric genus, i.e. $p_{g}(X)=0$, then the Albanese map \[
CH^{2}_{deg \ 0}(X) \to Alb(X)
\]
is an isomorphism, where $CH^{2}_{deg \ 0}(X)$ is the subgroup of $CH^{2}(X)$ consisting of zero cycles with degree zero and $Alb(X)$ is the Albanese variety.
\end{conjecture} 
This conjecture is closely related with the well-known example of Mumford (cf. Lewis \cite{Lewis1, Lewis2}, Mumford \cite{Mu}, Roitman \cite{Ro1, Ro2}, Voisin \cite{V1} et al) and had been studied intensively, for example, see Bloch-Kas-Lieberman \cite{BKL}, Bloch-Srinivas \cite{BS}, Hu \cite{Hu}, Pedrini-Weibel \cite{PW} and Voisin \cite{V2, V3}.

Stienstra \cite{Stien1} further studied these functors $\widetilde{CH}^{p}$ and $\widehat{CH}^{p}$, and computed $\widehat{CH}^{p}(A)$ in the case that $k$ is an extension of $\mathbb{Q}$ of finite transcendence degree. He also considered the parallel situation in positive characteristic and developed Cartier-Dieudonn{\'e} theory for Chow groups in \cite{Stien2}. An excellent summary of these results is given by Bloch \cite{Bl4} (Chapter 6).

Let $\overline{K}^{M}_{p}(O_{X_{A}})$ be the relative K-group, which is defined to be the kernel of the morphism $ K^{M}_{p}(O_{X_{A}}) \to K^{M}_{p}(O_{X})$. One challenge in studying $\widehat{CH}^{p}(A)$, which already appeared in the case of $p=2$, is computation of $\overline{K}^{M}_{p}(O_{X_{A}})$. To get a feeling of this, we recall that, for $A=k[\varepsilon]$ the ring of dual numbers, van der Kallen \cite{VanderKallen} computed that $\overline{K}^{M}_{2}(O_{X_{k[\varepsilon]}})$ is isomorphic to the sheaf of absolute K$\mathrm{\ddot{a}}$hler differentials $\Omega^{1}_{X/ \mathbb{Q}}$. Consequently, there is an isomorphism
\begin{equation}
\widehat{CH}^{2}(k[\varepsilon]) = H^{2}(X, \Omega^{1}_{X/ \mathbb{Q}}),
\end{equation}
and $\widehat{CH}^{2}(k[\varepsilon])$ is called the formal tangent space to Chow group $CH^{2}(X)$. Bloch \cite{Bl1-tan, Bl3} and Maazen-Stienstra \cite{MaazenStien} made further computations of relative K-groups. The following is one of basic results on understanding formal completion of $CH^{2}(X)$.
\begin{theorem}  [cf. Theorem 6.2 of \cite{Bl4}] \label{t:Bloch th}
Let $k$ be a field of characteristic zero. Let $R$ be a local $k$-algebra, and $A$ an augmented artinian $k$-algebra with augmentation ideal $m_{A}$. We write
$S=R \otimes_{k} A$ and $I=R \otimes_{k} m_{A}$, and define
\[
K_{2}(S,I)=\mathrm{kernel \ of \ }\{K_{2}(S) \to K_{2}(R) \},
\]
\[
\Omega^{1}_{S,I}=\mathrm{kernel \ of \ }\{\Omega^{1}_{S/ \mathbb{Q}} \to \Omega^{1}_{R/ \mathbb{Q}}\}.
\]
The universal derivation $d: S \to \Omega^{1}_{S/ \mathbb{Q}}$ induces $d: I \to \Omega^{1}_{S,I}$ and there is an isomorphism $K_{2}(S,I) \cong \Omega^{1}_{S,I}/dI$.

\end{theorem}

In the pioneering work \cite{GGtangentspace}, Green and Griffiths studied deformation of algebraic cycles of a smooth projective variety $X$ and investigated geometric meaning behind the formal tangent space to $CH^{2}(X)$ defined via the isomorphism (1.4). In particular, they computed the tangent space to zero cycles of a surface and justified that the formal tangent space to $CH^{2}(X)$ carried concrete geometric meaning, see Theorem 8.47 of \cite{GGtangentspace}. Inspired by a list of questions asked by Green and Griffiths in \cite{GGtangentspace}, Dribus, Hoffman and the author used higher K-theory to extend much of their theory in \cite{DHY, Y1, Y2, Y3,Y4}.
Especially relevant to the present paper is the following question in section 7.2 of \cite{GGtangentspace} (see also Question 1.2 in \cite{Y1}).
\begin{question}  [\cite{GGtangentspace}] \label{question: comparetangent}
Let $X$ be a smooth projective variety over a field $k$ of characteristic zero and let $Y \subset X$ be a closed subvariety of codimension $p$, is it possible to define a map from the tangent space $\mathrm{T}_{Y}\mathrm{\mathbb{H}ilb}^{p}(X)$ of the Hilbert scheme at the point $Y$ to the tangent space of the cycle group $\mathrm{T}Z^{p}(X)$
 \[
  \mathrm{T}_{Y}\mathrm{\mathbb{H}ilb}^{p}(X)  \to \mathrm{T}Z^{p}(X) ?
 \]
\end{question}
For $p=\mathrm{dim}(X)$, Green-Griffiths \cite{GGtangentspace} answered this question by studying deformations of zero cycles over the ring of dual numbers. Their method was generalized by the author \cite{Y1}.

The ring of dual numbers is a special local artinian $k$-algebra. Green and Griffiths' question inspires us to compare deformation of subvarieties with that of algebraic cycles (classes) over arbitrary local artinian $k$-algebras. Then we come to the following question suggested by Bloch in the introduction of \cite{Bl3} (page 406).
\begin{question}  [\cite{Bl3}]\label{q:Bloch's ques}
Let $X$ be a smooth projective variety over a field $k$ of characteristic zero and let $Y \subset X$ be a closed subvariety of codimension $p$, is there a natural transformation from local Hilbert functor $\mathrm{\mathbb{H}ilb}$ (recalled in Definition \ref{d:Hilbert} below) to the functor $\widetilde{CH}^{p}$ (see (1.2))
\[
\mathrm{\mathbb{H}ilb} \to \widetilde{CH}^{p}?
\]

\end{question}

This question is closely related with the following one suggested by Green-Griffiths on page 471 of \cite{GGChow}.
\begin{question}  [\cite{GGChow}]\label{q:GG's ques}
Let $X$ be a smooth projective variety over a field $k$ of characteristic zero. For $A \in Art_{k}$, we write $X_A=X \times _{\mathrm{Spec}(k)}\mathrm{Spec}(A)$. Is it possible to extend Bloch formula (1.1) from $X$ to its infinitesimal thickening $X_{A}$? In other words, do we have the following identification
\[
CH^{p}(X_{A})_{\mathbb{Q}}=H^{p}(X, K^{M}_{p}(O_{X_A}))_{\mathbb{Q}}?
\]
\end{question}
By modifying Balmer's tensor triangular Chow groups \cite{Ba}, we answered this question when $A$ is a truncated polynomial $k[t]/(t^{j})$ in \cite{Y2}.

Guided by Question \ref{q:Bloch's ques} and Question \ref{q:GG's ques}, this paper is organized as follows. In section 2, after recalling Bloch-Ogus theorem, cyclic homology and Milnor Chow groups, we answer Question \ref{q:GG's ques} in Theorem \ref{theorem:Agree Milnor Chow}. In the third section, we construct a map from local Hilbert functor to local homology in (3.2). With suitable assumptions, we use this map to answer Question \ref{q:Bloch's ques} in Theorem \ref{t: transf-Hilb-Chow}.

\textbf{Notation:}

(1). For any abelian group $M$, $M_{\mathbb{Q}}$ denotes $M \otimes_{\mathbb{Z}} \mathbb{Q}$. 

(2). If not stated otherwise, $k$ is a field of characteristic zero. Let $Sch/k$ be the category of schemes of finite type over $k$ and Let $Art_{k}$ denote the category of local artinian $k$-algebras with residue field $k$.

(3). If not stated otherwise, K-theory in this paper is Thomason-Trobaugh non-connective K-theory. For $X \in Sch/k$, let $Y \subset X$ be closed, Keller \cite{Keller, Keller2} defined cyclic homology complexes $HC(X)$ and $HC(X \ \mathrm{on} \  Y )$ from localization pairs (see Example 2.7 and 2.8 of \cite{CHSW} for details), which agree with the definitions of Weibel \cite{W3}. 

Following the convention in section 2 of \cite{CHSW}, we use cohomological notation for cyclic homology.

(4). For $F$ an abelian group-valued functor, we denote by $F(O_{X})$ the sheaf on a scheme $X \in Sch/k$ obtained by localizing $F$. The functor $F$ used in this paper are Milnor K-group $K^{M}_{*}(-)$, K-group $K_{*}(-)$, Hochschild homology $HH_{*}(-)$, cyclic homology $HC_{*}(-)$ and their eigenspaces of Adams operations $\psi^{m}$, denoted $K^{(l)}_{*}(-)$, $HH^{(l)}_{*}(-)$ and $HC^{(l)}_{*}(-)$ respectively.

\section{Bloch-Ogus theorem and local cohomology}

In this section, after recalling Bloch-Ogus theorem, cyclic homology and Milnor Chow group, we extend Bloch formula (1.1) in Theorem \ref{theorem:Agree Milnor Chow}.

\subsection{Bloch-Ogus Theorem}

Given a smooth algebraic variety $X$ and a cohomology theory $h$ satisfying natural axioms, the classical Bloch-Ogus theorem, says that the Zariski sheafification of the Cousin complex (formed from the coniveau spectral sequence) of $h$ is a flasque resolution of the Zariski sheaf associated to the presheaf $U \to h^{\ast}(U)$.

Bloch-Ogus \cite{BlochOgus} proved their theorem for \'etale cohomology with coefficients in roots of unity, by reducing to the ``effacement theorem" which was proved by using a geometric presentation lemma. Later, Gabber \cite{Gabber} gave a different proof of effacement theorem for \'etale cohomology.

In \cite{CTHK}, Colliot-Th\'el\`ene, Hoobler and Kahn axiomatized Gabber's proof and showed that his argument could be applied to any ``Cohomology theory with support'' which satisfies \'etale excision and a technical lemma (called ``Key lemma"). The latter follows either from homotopy invariance or from projective bundle formula. In particular, Gabber's argument works for K-theory and (negative) cyclic homology. This was used by Dribus, Hoffman and the author \cite{DHY} to study the deformation of algebraic cycles. We recall it briefly.

Let $h$ be a contravariant functor from the category $Sch/k$ to spectra or chain complexes. For $X \in Sch/k$, let $Y \subset X$ be closed, we can extend $h$ to the pair $(X,Y)$.
\begin{definition} \label{d: coh support}
For $h$ spectrum-valued, $h(X \ \mathrm{on} \ Y)$ is defined as the homotopy fiber of $h(X)\to h(X - Y)$. For any integer $p$, $h^{p}(X \ \mathrm{on} \ Y)$ is defined as homotopy group $\pi_{-p}(h(X \ \mathrm{on} \ Y))$.

For $h$ chain complex-valued, let $C_{\bullet}$ be the mapping cone of $h(X)\to h(X - Y)$, then $h(X \ \mathrm{on} \ Y)$ is defined as $C_{\bullet}[-1]$. For any integer $p$, $h^{p}(X \ \mathrm{on} \ Y)$ is defined as homology group $H_{-p}(h(X \ \mathrm{on} \ Y))$.
\end{definition}
This gives a ``cohomology theory with support" in the sense of Definition 5.1.1 of \cite{CTHK}. Following \cite{CTHK}, we recall \'etale excision and projective bundle formula.
\begin{definition} [\'Etale excision] \label{d:etale ex}
A functor $h$ is said to satisfy \'etale excision, if it is additive and if for any \'etale morphism $f: X' \to X$ such that $f^{-1}(Y) \xrightarrow{f} Y$ is an isomorphism with $Y \subset X$ closed, the pullback
\begin{equation*}
 f^{\ast} : h^{p}( X \ \mathrm{on} \ Y) \xrightarrow{\simeq} h^{p}( X^{'} \ \mathrm{on} \ f^{-1}(Y))
\end{equation*}
is an isomorphism for any integer $p$. 

The functor $h$ is said to satisfy Zariski excision 
if the pullback $f^{\ast}$ is an isomorphism for any integer $p$, when $f$ runs over all  open immersions.

\end{definition}

\begin{definition} [Projective bundle formula for projective line] \label{d:proj bundle formula}
The functor $h$ is said to satisfy projective bundle formula for projective line, if
\[
  h^{p}( X )\oplus h^{p}( X ) \xrightarrow{\simeq} h^{p}(\mathbb{P}_{X}^{1})
\]
is an isomorphism for any $X \in Sch/k$ and for any integer $p$, where $\mathbb{P}_{X}^{1}$ is the projective line over $X$.
\end{definition}

If the functor $h$ in Definition \ref{d: coh support} satisfies Zariski excision, then there exists a convergent spectral sequences, called coniveau spectral sequence (see section 1 of \cite{CTHK}),
\[
 E_{1}^{q,p} = \bigoplus_{x \in X^{(q)}}h^{q+p}(X \ \mathrm{on} \ x) \Longrightarrow h^{q+p}(X),
\]
where $X^{(q)}$ denotes the set of points of codimension $q$ in $X$ and
\[
h^{q+p}(X \ \mathrm{on} \ x) = \varinjlim_{x \in U}h^{q+p}(U \ \mathrm{on} \ \overline{\{x\}}\cap U).
\]
The $E_{1}$-terms give rise to Cousin complex of $h$
{\small
\begin{equation}
 0 \to \bigoplus_{x \in X^{(0)}}h^{p}(X \ \mathrm{on} \ x) \to \bigoplus_{x \in X^{(1)}}h^{p+1}(X \ \mathrm{on} \ x) \to \cdots.
\end{equation}
}

￼￼
The following setting is used below.
\begin{setting} \label{s:setting}
Let $X$ be a $d$-dimensional smooth projective variety over a field $k$ of characteristic zero, with generic point $\eta$. 

For $A \in Art_{k}$, we write $X_A=X \times _{\mathrm{Spec}(k)}\mathrm{Spec}(A)$. Let $F$ be a functor as in Notation (4) on page 4, we denote by $\overline{F}(O_{X_A})$ the kernel of the morphism (induced by augmentation $A \to k$) $F(O_{X_A}) \to F(O_{X})$. 
\end{setting}

\begin{theorem} [Bloch-Ogus Theorem] \label{t:Bloch-Ogus}
In notation of Setting \ref{s:setting}, if a functor $h$ in Definition \ref{d: coh support} satisfies \'etale excision and projective bundle formula for projective line, then for any integer $p$, the Zariski sheafification of the Cousin complex (2.1) is a flasque resolution of the sheaf associated to the presheaf $U \to h^{p}(O_{X}(U))$.
\end{theorem}

\begin{proof}
This was originally proved by Bloch-Ogus \cite{BlochOgus} for \'etale cohomology and it was extended by Gabber \cite{Gabber}. Colliot-Th\'el\`ene, Hoobler and Kahn applied Gabber's method to prove the theorem in a general context, see Corollary 5.1.11 and Proposition 5.4.3 of \cite{CTHK}.

\end{proof}

Universal exactness was originally introduced by Grayson \cite{Gray}. For arbitrary scheme $T \in Sch/k$ ($T$ might be singular), we can derive a new functor $h^{T}$ from the functor $h$ in Definition \ref{d: coh support},
\[
 h^{T}: \ X \to h(X \times_{\mathrm{Spec}(k)} T).
\]
If the functor $h$ satisfies \'etale excision and projective bundle formula for projective line, so does the new functor $h^{T}$. This implies that,
\begin{corollary} [Corollary 6.2.4 of \cite{CTHK}] \label{c:univ exact}
In notation of Setting \ref{s:setting}, if a functor $h$ in Definition \ref{d: coh support} satisfies \'etale excision and projective bundle formula for projective line, then for any integer $p$, the Zariski sheafification of the Cousin complex (2.1) of $h^{T}$ is a flasque resolution of the sheaf associated to the presheaf $U \to h^{p}(U \times_{\mathrm{Spec}(k)} T)$.
\end{corollary}

\begin{lemma} \label{l: K and HC effac}
Both K-theory and cyclic homology satisfy \'etale excision and projective bundle formula.
\end{lemma}

\begin{proof}
For K-theory, it was proved in Theorem 7.1 (for \'etale excision)  and in Theorem 7.3 (for projective bundle formula) of \cite{TT}. For cyclic homology, see Example 2.8 of \cite{CHSW} (for \'etale excision) and Remark 2.11 of \cite{CHSW} (for projective bundle formula).

\end{proof}

When the functor $h$ is the K-theory spectrum $\mathcal{K}(X)$, the associated Cousin complex (2.1) is the Bloch-Gersten-Quillen sequence
{\footnotesize
\begin{align}
0  \rightarrow  \bigoplus_{x \in X^{(0)}}K_{p}(O_{X,x}) \to \bigoplus_{x \in X^{(1)}}K_{p-1}(O_{X,x} \ \mathrm{on} \ x)  \to \cdots.
\end{align}
}

\begin{corollary} \label{c:BGQresolution}
In Setting \ref{s:setting}, for each integer $p \geq 0$, the Zariski sheafifications of the Bloch-Gersten-Quillen sequences (2.2) of $X$ and $X_A$ are flasque resolutions of the K-theory sheaves $K_{p}(O_X)$ and $K_{p}(O_{X_A})$ respectively.
\end{corollary}

\begin{proof}
By Quillen's d\'evissage, the Bloch-Gersten-Quillen sequence (2.2) of $X$ has the form
{\footnotesize
 \begin{align*}
0  \rightarrow K_{p}(O_{X,\eta}) \rightarrow \dots \rightarrow \bigoplus_{x \in X^{(p-1)}}K_{1}(k(x))  \to \bigoplus_{x \in X^{(p)}}K_{0}(k(x))  \to 0,
\end{align*}
}whose Zariski sheafification is a flasque resolution of the sheaf $K_{p}(O_X)$. This was used by Quillen \cite{Quillen} to prove Bloch formula.

By Corollary \ref{c:univ exact} and Lemma \ref{l: K and HC effac}, the Zariski sheafification of the Bloch-Gersten-Quillen sequence (2.2) of $X_A$, which has the form
{\small
\begin{align}
& 0  \to K_{p}(O_{X_A,\eta}) \rightarrow \dots \rightarrow \bigoplus_{x \in X^{(p)}}K_{0}(O_{X_A,x} \ \mathrm{on} \ x)  \\
& \to \bigoplus_{x \in X^{(p+1)}}K_{-1}(O_{X_A,x} \ \mathrm{on} \ x)   \rightarrow \cdots,  \notag
\end{align}
}is a flasque resolution of the sheaf $K_{p}(O_{X_{A}})$.

\end{proof}

It is worth noting that nontrivial negative K-groups may appear in the sequence (2.3).

The closed immersion $X \to X_{A}$ induces a map between K-theory spectra $\mathcal{K}(X_{A}) \to \mathcal{K}(X)$, which further induces a map between Bloch-Gersten-Quillen sequences. Since the closed immersion $X \to X_{A}$ has a section $X_{A} \to X$, there is a split commutative diagram
{\tiny
\begin{equation}
  \begin{CD}
     0 @. 0 @. 0\\
      @VVV @VVV @VVV\\
     \overline{K}_{p}(O_{X_{A},\eta}) @<<< K_{p}(O_{X_{A},\eta}) @>>>  K_{p}(O_{X,\eta}) \\
     @VVV @VVV @VVV\\
      \bigoplus\limits_{x \in X ^{(1)}} \overline{K}_{p-1}(O_{X_{A},x} \ \mathrm{on} \ x) @<<< \bigoplus\limits_{x \in X^{(1)}}K_{p-1}(O_{X_{A},x} \ \mathrm{on} \ x) @>>>  \bigoplus\limits_{x \in X ^{(1)}}K_{p-1}(O_{X,x} \ \mathrm{on} \ x) \\
    @VVV  @VVV @VVV\\
     \vdots @<<<  \vdots @>>> \vdots \\ 
      @VVV @VVV @VVV\\
     \bigoplus\limits_{x \in X ^{(d)}} \overline{K}_{p-d}(O_{X_{A},x} \ \mathrm{on} \ x) @<<<  \bigoplus\limits_{x \in X^{(d)}}K_{p-d}(O_{X_{A},x} \ \mathrm{on} \ x) @>>> \bigoplus\limits_{x \in X ^{(d)}}K_{p-d}(O_{X,x} \ \mathrm{on} \ x) \\
     @VVV @VVV @VVV\\
     0 @. 0 @. 0,
  \end{CD}
\end{equation}
}where each $\overline{K}_{*}(O_{X_{A},x} \ \mathrm{on} \ x)$ is the kernel of the map (induced by augmentation $A \to k$) $K_{*}(O_{X_{A},x} \ \mathrm{on} \ x) \to K_{*}(O_{X,x} \ \mathrm{on} \ x)$. Let $\overline{\mathcal{K}}(X_{A})$ denote the homotopy fiber of $\mathcal{K}(X_{A}) \to \mathcal{K}(X)$, the left column of the diagram (2.4) is the Cousin complex of the spectrum $\overline{\mathcal{K}}(X_{A})$.

\subsection{Cyclic homology}
Hochschild and cyclic homology are defined over $\mathbb{Q}$ here. For $R$ a commutative $\mathbb{Q}$-algebra, Hochschild homology $HH_{*}(R)$ and cyclic homology $HC_{*}(R)$ carry Lambda operations $\lambda^{m}$ and Adams operations $\psi^{m}$, see section 4.5 of \cite{Loday} and section 9.4.3 of \cite{W2} for details. In fact, the action of symmetric group naturally splits Hochschild complex $HH(R)$ and cyclic homology complex $HC(R)$ into sums of sub-complexes $HH^{(l)}(R)$ and $HC^{(l)}(R)$ respectively. For each integer $p \geq 1$, this decomposes $HH_{p}(R)$ and $HC_{p}(R)$ into direct sums of eigenspaces
\begin{equation*}
HH_{p}(R)=HH^{(1)}_{p}(R) \oplus \cdots \oplus HH^{(p)}_{p}(R),
\end{equation*}
\begin{equation}
HC_{p}(R)=HC^{(1)}_{p}(R) \oplus \cdots \oplus HC^{(p)}_{p}(R).
\end{equation}
For $p=0$, $HC_{0}(R)=HC^{(0)}_{0}(R)=R$.

\begin{lemma} [Ex 9.4.4 and Corollary 9.8.16 of \cite{W2}]\label{lemma: l-l-omega}
With notation as above, there are isomorphisms
\[
HH^{(p)}_{p}(R) \cong \Omega^{p}_{R/\mathbb{Q}}, \ HC^{(p)}_{p}(R)=\dfrac{\Omega^{p}_{R/\mathbb{Q}}}{d\Omega^{p-1}_{R/\mathbb{Q}}}.
\]

\end{lemma}

These operations $\lambda^{m}$ and $\psi^{m}$ can be extended to cyclic homology $HC_{*}(X)$, where $X \in Sch/k$, see Weibel \cite{W4}. Let $Y \subset X$ be closed, since cyclic homology satisfies Zariski descent, we can identify cyclic homology $HC_{*}(X \ \mathrm{on} \ Y)$ with hypercohomology
\[
HC_{*}(X \ \mathrm{on} \ Y) = \mathbb{H}^{-*}_{Y}(X, HC(X)),
\] 
where $HC(X)$ is the cyclic homology complex of $X$.  This enables us to further extend $\lambda^{m}$ and $\psi^{m}$ to $HC_{*}(X \ \mathrm{on} \ Y)$.

Combing Corollary \ref{c:univ exact} with Lemma \ref{l: K and HC effac}, one has
\begin{lemma}  \label{l:HC-flasque sheaf}
In Setting \ref{s:setting}, for each integer $p \geq 0$, the Zariski sheafification of the following Cousin complex of cyclic homology of $X_{A}$
{\small
\begin{align}
0  \to HC_{p}(O_{X_A,\eta}) \to \bigoplus_{x \in X^{(1)}}HC_{p-1}(O_{X_A,x} \ \mathrm{on} \ x) \to \cdots,
\end{align}
}is a flasque resolution of the sheaf $HC_{p}(O_{X_{A}})$.

\end{lemma}

The differentials of the complex of (2.6) respect Adams operations $\psi^{m}$. This yields that

\begin{lemma} \label{l: Eigen-HC-flasque}
In Setting \ref{s:setting}, for each integer $p \geq 0$, the Zariski sheafification of the complex
{\small
\begin{align}
0  \to  HC^{(l)}_{p}(O_{X_A,\eta}) \to \bigoplus_{x \in X^{(1)}}HC^{(l)}_{p-1}(O_{X_A,x} \ \mathrm{on} \ x) \to \cdots,
\end{align}
}is a flasque resolution of the sheaf $HC^{(l)}_{p}(O_{X_{A}})$, where the integer $l$ satisfying that $0 \leq l \leq p$ and each $HC^{(l)}_{*}(-)$ is eigenspace of Adams operations $\psi^{m}$.
 
\end{lemma}

We are mainly interested in the case $l=p$ below. Let $q$ be an integer satisfying that $1 \leq q \leq d$, where $d =\mathrm{dim}(X)$. For $x \in X^{(q)}$, let $\overline{HC}^{(p)}_{p-q}(O_{X_{A},x} \ \mathrm{on} \ x)$ be the kernel of the map (induced by $A \to k$) $HC^{(p)}_{p-q}(O_{X_{A},x} \ \mathrm{on} \ x) \to HC^{(p)}_{p-q}(O_{X,x} \ \mathrm{on} \ x)$. Let $l=p$ in (2.7), Lemma \ref{l: Eigen-HC-flasque} implies that 

\begin{corollary} \label{c:overline-HC-p}
In Setting \ref{s:setting}, for each integer $p \geq 0$, the Zariski sheafification of the complex
{\small
\begin{align}
0  \to \overline{HC}^{(p)}_{p}(O_{X_A,\eta}) \to \bigoplus_{x \in X^{(1)}}\overline{HC}^{(p)}_{p-1}(O_{X_A,x} \ \mathrm{on} \ x) \to \cdots,
\end{align}
}is a flasque resolution of the sheaf $\overline{HC}^{(p)}_{p}(O_{X_{A}})$.
\end{corollary}
 We want to compute each group $\overline{HC}^{(p)}_{*}(O_{X_{A},x} \ \mathrm{on} \ x)$ of the complex (2.8). We refer to chapter IV of \cite{Hart1} for definitions and properties of local cohomologies of abelian sheaves.

\begin{lemma} \label{l: compute-HC-relative}
In Setting \ref{s:setting}, let $q$ be an integer satisfying that $1 \leq q \leq d$, where $d =\mathrm{dim}(X)$. For $x \in X^{(q)}$ and for each integer $p \geq 0$, {\small $\overline{HC}^{(p)}_{p-q}(O_{X_{A},x} \ \mathrm{on} \ x)$}\footnote{The index $p-q$ might be negative.} is isomorphic to local cohomology $H^{q}_{x}(\overline{HC}^{(p)}_{p}(O_{X_{A}}))$
\[
\overline{HC}^{(p)}_{p-q}(O_{X_{A},x} \ \mathrm{on} \ x)=H^{q}_{x}(\overline{HC}^{(p)}_{p}(O_{X_{A}})).
\]
\end{lemma}
 As recalled in the beginning of section 2.2, the cyclic homology complexes $HC(O_{X_{A},x})$ and $HC(O_{X,x})$ split into direct sums of sub-complexes $HC^{(l)}(O_{X_{A},x})$ and $HC^{(l)}(O_{X,x})$ respectively. We are interested in $HC^{(p)}(O_{X_{A},x})$ and $HC^{(p)}(O_{X,x})$, and denote by $\overline{HC}^{(p)}(O_{X_{A},x})$ the kernel of the natural map of complexes
 \[
 HC^{(p)}(O_{X_{A},x}) \to HC^{(p)}(O_{X,x}).
 \]

\begin{proof}
Since cyclic homology satisfies Zariski descent, we can identify $\overline{HC}^{(p)}_{p-q}(O_{X_{A},x} \ \mathrm{on} \ x)$ with hypercohomogy
\[
\overline{HC}^{(p)}_{p-q}(O_{X_{A},x} \ \mathrm{on} \ x)= \mathbb{H}^{-(p-q)}_{x}(O_{X,x}, \overline{HC}^{(p)}(O_{X_{A},x})).
\]

There exists a spectral sequence
{\tiny
\[
E_{2}^{i, j} = H^{i}_{x}(O_{X,x}, H^{j}(\overline{HC}^{(p)}(O_{X_{A},x})))
 \Longrightarrow \mathbb{H}^{-(p-q)}_{x}(O_{X,x}, \overline{HC}^{(p)}(O_{X_{A},x})),
\]
}where $i+j=-(p-q)$. We use cohomological notation for cyclic homology (see Notation (3) on page 4), so $H^{j}(\overline{HC}^{(p)}(O_{X_{A},x}))=\overline{HC}^{(p)}_{-j}(O_{X_{A},x})$, the above spectral sequence can be rewritten as 
{\tiny
\begin{equation}
E_{2}^{i, j} = H^{i}_{x}(O_{X,x}, \overline{HC}^{(p)}_{j}(O_{X_{A},x}))
 \Longrightarrow \mathbb{H}^{-(p-q)}_{x}(O_{X,x}, \overline{HC}^{(p)}(O_{X_{A},x})),
\end{equation}
}where $i$ and $j$ are non-negative integers, and $i-j=-(p-q)$.

Since the Krull dimension of $O_{X,x}$ is $q$, if $i > q$, then the local cohomology $H^{i}_{x}(O_{X,x}, \overline{HC}^{(p)}_{j}(O_{X_{A},x}))=0$ for each $j$. This shows that the index $i$ in non-zero terms of the spectral sequence (2.9) satisfies that $0 \leqslant i \leqslant q$. It follows that $j=i+p-q \leqslant p$. If $j < p$, then $\overline{HC}^{(p)}_{j}(O_{X_{A},x})=0$, see (2.5) on page 8. Hence, the index $j=i+p-q$ in non-zero terms of the spectral sequence (2.9) can only be $p$, which implies that $i=q$.

In conclusion, the only non-zero term in the spectral sequence (2.9) is $H^{q}_{x}(O_{X,x}, \overline{HC}^{(p)}_{p}(O_{X_{A},x}))$. So the spectral sequence (2.9) degenerates and
\[
\overline{HC}^{(p)}_{p-q}(O_{X_{A},x} \ \mathrm{on} \ x)=H^{q}_{y}(O_{X,x},\overline{HC}^{(p)}_{p}(O_{X_{A},x})).
\]

\end{proof}  

\begin{definition} [cf. Def 3.2 of \cite{BEK2}]\label{d: CM}
 For $X \in Sch/k$, an abelian sheaf $F$ on $X$ is called Cohen-Macaulay, if for every scheme point $y \in X$, $H_{y}^{i}(X, F)=0$ for $i \neq \mathrm{codim\{y\}}$.
 \end{definition}

To see the importance of Cohen-Macaulay sheaves, we recall that, for $X \in Sch/k$ and for an abelian sheaf $F$ on $X$, the Cousin complex of $F$ constructed in Proposition 2.3 of chapter IV of \cite{Hart1} has the form
\begin{equation}
0 \to \bigoplus_{x \in X^{(0)}}H^{0}_{x}(F) \to \bigoplus_{x \in X^{(1)}}H^{1}_{x}(F) \to \cdots.
\end{equation}

\begin{lemma} [Prop 2.6 of chapter IV of \cite{Hart1}]\label{l:CM-Resolution}
The following are equivalent:
\begin{itemize}
\item [$\mathrm{(1)}$] the abelian sheaf $F$ is Cohen-Macaulay,
\item [$\mathrm{(2)}$] the Zariski sheafification of the Cousin complex (2.10) of $F$ is a flasque resolution of $F$.
\end{itemize}

\end{lemma}

For $X$ a smooth projective variety over a field $k$ of characteristic zero, $\Omega^{*}_{X/k}$ is Cohen-Macaulay (see page 239 of \cite{Hart1}). Kerz \cite{Kerz} proved that Milnor K-theory sheaf $K^{M}_{*}(O_{X})$ is Cohen-Macaulay. Bloch-Esnault-Kerz generalized these examples in Prop.3.5 of \cite{BEK2} and applied it to the infinitesimal  study of Chow groups. This motivates us to find more examples of Cohen-Macaulay sheaves.

In Setting \ref{s:setting}, by Corollary \ref{c:overline-HC-p} and Lemma \ref{l: compute-HC-relative}, the Zariski sheafification of the Cousin complex (2.10) of the sheaf $\overline{HC}^{(p)}_{p}(O_{X_{A}})$
\begin{equation*}
0 \to \overline{HC}^{(p)}_{p}(O_{X_{A}, \eta}) \to \bigoplus_{x \in X^{(1)}}H^{1}_{x}(\overline{HC}^{(p)}_{p}(O_{X_{A}})) \to \cdots,
\end{equation*}is a flasque resolution of $\overline{HC}^{(p)}_{p}(O_{X_{A}})$. It follows from Lemma \ref{l:CM-Resolution} that 

\begin{corollary} \label{c: HC-CM}
In Setting \ref{s:setting}, for each integer $p \geq 0$, the sheaf $\overline{HC}^{(p)}_{p}(O_{X_{A}})$ is Cohen-Macaulay.

\end{corollary}

\begin{lemma} \label{l: HH-CM}
In Setting \ref{s:setting}, for each integer $p \geq 0$, the sheave $\overline{HH}^{(p)}_{p}(O_{X_{A}})$ is Cohen-Macaulay.
\end{lemma}

\begin{proof}
By Lemma \ref{lemma: l-l-omega}, $HH^{(p)}_{p}(O_{X_{A}}) \cong \Omega^{p}_{X_A/ \mathbb{Q}}$. We first show that $HH^{(p)}_{p}(O_{X_{A}})$ is Cohen-Macaulay. Let $F^{j}$ be the image of the map $\Omega^{j}_{k/ \mathbb{Q}} \otimes_{k} \Omega^{p-j}_{X_A/ \mathbb{Q}} \to \Omega^{p}_{X_A/ \mathbb{Q}}$, where $j=0, 1, \cdots, p$. There is a filtration on $\Omega^{p}_{X_A/ \mathbb{Q}}$ given by
\begin{equation*}
\Omega^{p}_{X_A/ \mathbb{Q}}=F^{0} \supset F^{1} \supset \cdots \supset F^{p} \supset F^{p+1}=0,
\end{equation*}
whose associated graded piece is $Gr^{j}\Omega^{p}_{X_A/ \mathbb{Q}}=F^{j}/F^{j+1}=\Omega^{j}_{k/ \mathbb{Q}} \otimes_{k} \Omega^{p-j}_{X_A/ k}$. In particular, $Gr^{p}\Omega^{p}_{X_A/ \mathbb{Q}}=F^{p}=\Omega^{p}_{k/ \mathbb{Q}} \otimes_{k} O_{X_A}$.

There is an isomorphism of sheaves
\[
\Omega^{p-j}_{X_A/ k}= \bigoplus_{j_{1}+j_{2}=p-j}\Omega^{j_{1}}_{X/ k} \otimes_{k}\Omega^{j_{2}}_{A/ k},
\]
which can be checked locally ($X_A$ and $X$ have the same underlying space). Since each $\Omega^{j_{1}}_{X/ k}$ is Cohen-Macaulay, so is $\Omega^{p-j}_{X_A/ k}$. It follows that each $Gr^{j}\Omega^{p}_{X_A/ \mathbb{Q}}=\Omega^{j}_{k/ \mathbb{Q}} \otimes_{k} \Omega^{p-j}_{X_A/ k}$ is Cohen-Macaulay.

There is a short exact sequence
\[
0 \to F^{p} \to F^{p-1} \to Gr^{p-1} \to 0,
\]
where both $F^{p}=Gr^{p}\Omega^{p}_{X_A/ \mathbb{Q}}=\Omega^{p}_{k/ \mathbb{Q}} \otimes_{k} O_{X_A}$ and $Gr^{p-1}=\Omega^{p-1}_{k/ \mathbb{Q}} \otimes_{k} \Omega^{1}_{X_A/ k}$ are Cohen-Macaulay, so the associated long exact sequence of local cohomology implies that $F^{p-1}$ is Cohen-Macaulay. We are able to prove that each $F^{j}$  is Cohen-Macaulay by continuing this procedure. In particular, $F^{0}=\Omega^{p}_{X_A/ \mathbb{Q}}$ is Cohen-Macaulay. When $A=k$, $HH^{(p)}_{p}(O_{X}) \cong \Omega^{p}_{X/ \mathbb{Q}}$ is Cohen-Macaulay.

The short exact sequence
\[
 0 \to \overline{HH}^{(p)}_{p}(O_{X_{A}}) \to HH^{(p)}_{p}(O_{X_{A}}) \to HH^{(p)}_{p}(O_{X}) \to 0
\]
is split, where both $HH^{(p)}_{p}(O_{X_{A}})$ and $HH^{(p)}_{p}(O_{X})$ are Cohen-Macaulay, so is $\overline{HH}^{(p)}_{p}(O_{X_{A}})$.
\end{proof}

\subsection{Milnor Chow groups}

For $X \in Sch/k$, it is well known that Grothendieck group of $X$ carries Adams operations $\psi^{m}$, which is induced from exterior powers of vector bundles on $X$. These operations can be extended to higher algebraic K-theory.
For $Y \subset X$ a closed subscheme, Soul\'e \cite{Soule} and Levine \cite{Levine} defined Adams operations $\psi^{m}$ on K-groups with supports $K_{n}( X \ \mathrm{on} \ Y)$, where $n \geq 0$. 

Since the appearance of nontrivial negative K-groups in our study, we need to extend Adams operations $\psi^{m}$ to negative range. According to Weibel \cite{W1} (section 8), this can be done inductively by using
Bass fundamental exact sequence. We have used this method in \cite{DHY} (section 8.2).

Let $I$ be a nilpotent ideal in a commutative $\mathbb{Q}$-algebra $R$. We define the relative K-group $K_{n}(R,I)$ to be the kernel of the morphism $K_{n}(R) \to K_{n}(R/I)$ and define $K_{n}^{(l)}(R,I)$ to be the eigenspace of $\psi^{m}=m^{l}$, where $\psi^{m}$ is Adams operations on $K_{n}(R,I)$. The relative cyclic homology $HC_{n-1}(R,I)$ and $HC_{n-1}^{(l-1)}(R,I)$ are defined similarly. Goodwillie and Cathelineau proved that these relative groups are connected by the relative Chern character.
\begin{theorem} [\cite{Cath,Good}] \label{theorem: GoodwillieCathelineau}
With notation as above, the relative Chern character induces an isomorphism between $K_{n}(R,I)_{\mathbb{Q}}$ and $HC_{n-1}(R,I)$, which respects Adams operations 
\[
K_{n}(R,I)_{\mathbb{Q}} \xrightarrow{\cong} HC_{n-1}(R,I), \ K_{n}^{(l)}(R,I)_{\mathbb{Q}}  \xrightarrow{\cong} HC_{n-1}^{(l-1)}(R,I).
\]
\end{theorem}

This theorem is very useful to compute relative K-groups. Corti$\mathrm{\tilde{n}}$as-Haesemeyer-Weibel \cite{CHW} generalized it to space level. We adopt it to Setting \ref{s:setting} and refer to appendix B of \cite{CHW} for a general form. 

Let $\mathcal{HC}(X_{A})$ and $\mathcal{HC}(X)$ be the Eilenberg-Mac Lane spectra associated to cyclic homogy complexes $HC(X_{A})$ and $HC(X)$ respectively. Let $\overline{\mathcal{HC}}(X_{A})$ be the homotopy fiber of $\mathcal{HC}(X_{A}) \to \mathcal{HC}(X)$, we define $\overline{\mathcal{HC}}^{(l-1)}(X_{A})$ as the homotopy fiber of the map $\psi^{m}-m^{l}$ on $\overline{\mathcal{HC}}(X_{A})$. The spectra $\overline{\mathcal{K}}(X_{A})$ and $\overline{\mathcal{K}}^{(l)}(X_{A})$ are defined similarly. 

Theorem \ref{theorem: GoodwillieCathelineau} can be generalized in the following way.
\begin{theorem} [cf. Theorem B.11 of \cite{CHW}] \label{theorem: CHW}
In Setting \ref{s:setting}, the relative Chern character induces homotopy equivalence of spectra
 \[
 \overline{\mathcal{K}}(X_{A}) \xrightarrow{\simeq} \overline{\mathcal{HC}}(X_{A})[1], \ \overline{\mathcal{K}}^{(l)}(X_{A})  \xrightarrow{\simeq}  \overline{\mathcal{HC}}^{(l-1)}(X_{A})[1].
 \]
\end{theorem}

For each integer $p \geq 1$, the homotopy equivalence of spectra 
\[
\overline{\mathcal{K}}^{(l)}(X_{A}) \xrightarrow{\simeq} \overline{\mathcal{HC}}^{(l-1)}(X_{A})[1]
\]
induces a commutative diagram of Cousin complexes
{\footnotesize
\begin{equation}
  \begin{CD}
     0 @. 0 \\
      @VVV @VVV \\
     \overline{K}^{(l)}_{p}(O_{X_{A},\eta})_{\mathbb{Q}} @>\cong>> \overline{HC}^{(l-1)}_{p-1}(O_{X_{A},\eta}) \\
     @VVV @VVV \\
      \bigoplus\limits_{x \in X ^{(1)}} \overline{K}^{(l)}_{p-1}(O_{X_{A},x} \ \mathrm{on} \ x)_{\mathbb{Q}} @>\cong>> \bigoplus\limits_{x \in X ^{(1)}} \overline{HC}^{(l-1)}_{p-2}(O_{X_{A},x} \ \mathrm{on} \ x)\\
    @VVV  @VVV \\
     \vdots @>>>  \vdots \\ 
      @VVV @VVV \\
     \bigoplus\limits_{x \in X ^{(d)}} \overline{K}^{(l)}_{p-d}(O_{X_{A},x} \ \mathrm{on} \ x)_{\mathbb{Q}} @>\cong>>  \bigoplus\limits_{x \in X ^{(d)}} \overline{HC}^{(l-1)}_{p-d-1}(O_{X_{A},x} \ \mathrm{on} \ x)\\
     @VVV @VVV \\
     0 @. 0.
  \end{CD}
\end{equation}
}We explain the notations of diagram (2.11) briefly. Let $K^{(l)}_{*}(O_{X_{A},x} \ \mathrm{on} \ x)$ and $K^{(l)}_{*}(O_{X,x} \ \mathrm{on} \ x)$ denote eigenspaces of Adams operations $\psi^{m}=m^{l}$ respectively, $\overline{K}^{(l)}_{*}(O_{X_{A},x} \ \mathrm{on} \ x)$ is the kernel of the morphism 
\[
K^{(l)}_{*}(O_{X_{A},x} \ \mathrm{on} \ x) \to K^{(l)}_{*}(O_{X,x} \ \mathrm{on} \ x),
\]and $\overline{HC}^{(l-1)}_{*}(O_{X_{A},x} \ \mathrm{on} \ x)$ is defined similarly.

Combing Corollary \ref{c:BGQresolution}, Lemma \ref{l: Eigen-HC-flasque} with diagrams (2.4) and (2.11), one has
\begin{lemma}  \label{l:eigen-i-diag}
In Setting \ref{s:setting}, there is a commutative diagram
{\tiny
\[
  \begin{CD}
     0 @. 0 @. 0\\
      @VVV @VVV @VVV\\
     \overline{HC}^{(l-1)}_{p-1}(O_{X_{A},\eta}) @<<< K^{(l)}_{p}(O_{X_{A},\eta})_{\mathbb{Q}} @>>>  K^{(l)}_{p}(O_{X,\eta})_{\mathbb{Q}}  \\
     @VVV @VVV @VVV\\
      \bigoplus\limits_{x \in X ^{(1)}} \overline{HC}^{(l-1)}_{p-2}(O_{X_{A},x} \ \mathrm{on} \ x) @<<< \bigoplus\limits_{x \in X^{(1)}}K^{(l)}_{p-1}(O_{X_{A},x} \ \mathrm{on} \ x)_{\mathbb{Q}}  @>>>  \bigoplus\limits_{x \in X ^{(1)}}K^{(l)}_{p-1}(O_{X,x} \ \mathrm{on} \ x)_{\mathbb{Q}} \\
    @VVV  @VVV @VVV\\
     \vdots @<<<  \vdots @>>> \vdots \\ 
      @VVV @VVV @VVV\\
     \bigoplus\limits_{x \in X ^{(d)}} \overline{HC}^{(l-1)}_{p-d-1}(O_{X_{A},x} \ \mathrm{on} \ x) @<<<  \bigoplus\limits_{x \in X^{(d)}}K^{(l)}_{p-d}(O_{X_{A},x} \ \mathrm{on} \ x)_{\mathbb{Q}} @>>> \bigoplus\limits_{x \in X ^{(d)}}K^{(l)}_{p-d}(O_{X,x} \ \mathrm{on} \ x)_{\mathbb{Q}} \\
     @VVV @VVV @VVV\\
     0 @. 0 @. 0,
  \end{CD}
\]
}where each column is a complex whose Zariski sheafification is a flasque resolution of $\overline{HC}^{(l-1)}_{p-1}(O_{X_A})$, $K^{(l)}_{p}(O_{X_A})_{\mathbb{Q}}$ and $K^{(l)}_{p}(O_{X})_{\mathbb{Q}}$ respectively.

\end{lemma}

Combining Lemma \ref{l:eigen-i-diag} (let $l=p$) with Lemma \ref{l: compute-HC-relative}, one has
\begin{theorem} \label{theorem: Main-Diag}
In Setting \ref{s:setting}, for each integer $p \geq 1$, there exists the following commutative diagram
{\tiny
\[
  \begin{CD}
     0 @. 0 @. 0\\
      @VVV @VVV @VVV\\
      \overline{HC}^{(p-1)}_{p-1}(O_{X_{A},\eta}) @<\mathrm{Ch}<< K^{(p)}_{p}(O_{X_{A},\eta})_{\mathbb{Q}} @>\mathrm{Pr}>>  K^{(p)}_{p}(O_{X,\eta})_{\mathbb{Q}} \\
     @VVV @VVV @VVV\\
      \bigoplus\limits_{x \in X ^{(1)}} H^{1}_{x}(\overline{HC}^{(p-1)}_{p-1}(O_{X_{A}})) @<<< \bigoplus\limits_{x \in X^{(1)}}K^{(p)}_{p-1}(O_{X_{A},x} \ \mathrm{on} \ x)_{\mathbb{Q}} @>>>  \bigoplus\limits_{x \in X ^{(1)}}K^{(p)}_{p-1}(O_{X,x} \ \mathrm{on} \ x)_{\mathbb{Q}} \\
    @VVV  @VVV @VVV\\
     \vdots @<<<  \vdots @>>> \vdots \\ 
      @VVV @VVV @VVV\\
     \bigoplus\limits_{x \in X ^{(p-1)}} H^{p-1}_{x}(\overline{HC}^{(p-1)}_{p-1}(O_{X_{A}})) @<<<  \bigoplus\limits_{x \in X^{(p-1)}}K^{(p)}_{1}(O_{X_{A},x} \ \mathrm{on} \ x)_{\mathbb{Q}}
      @>>> \bigoplus\limits_{x \in X^{(p-1)}}K^{(p)}_{1}(O_{X,x} \ \mathrm{on} \ x)_{\mathbb{Q}} \\
     @V \partial_{1,X_{A}}^{p-1,-p}VV @Vd_{1,X_{A}}^{p-1,-p}VV @Vd_{1,X}^{p-1,-p}VV\\
      \bigoplus\limits_{x \in X ^{(p)}} H^{p}_{x}(\overline{HC}^{(p-1)}_{p-1}(O_{X_{A}})) @<<<  \bigoplus\limits_{x \in X^{(p)}}K^{(p)}_{0}(O_{X_{A},x} \ \mathrm{on} \ x)_{\mathbb{Q}}
      @>>> \bigoplus\limits_{x \in X ^{(p)}}K^{(p)}_{0}(O_{X,x} \ \mathrm{on} \ x)_{\mathbb{Q}} \\
     @V \partial_{1,X_{A}}^{p,-p}VV  @Vd_{1,X_{A}}^{p,-p}VV @Vd_{1,X}^{p,-p}VV\\
     \bigoplus\limits_{x \in X ^{(p+1)}} H^{p+1}_{x}(\overline{HC}^{(p-1)}_{p-1}(O_{X_{A}})) @<<< \bigoplus\limits_{x \in X^{(p+1)}}K^{(p)}_{-1}(O_{X_{A},x} \ \mathrm{on} \ x)_{\mathbb{Q}}
      @>>> \bigoplus\limits_{x \in X ^{(p+1)}}K^{(p)}_{-1}(O_{X,x} \ \mathrm{on} \ x)_{\mathbb{Q}} \\
      @VVV  @VVV @VVV\\
     \vdots @<<< \vdots @>>> \vdots \\ 
     @VVV @VVV @VVV\\
     \bigoplus\limits_{x \in X ^{(d)}} H^{d}_{x}(\overline{HC}^{(p-1)}_{p-1}(O_{X_{A}}))@<<<  \bigoplus\limits_{x \in X^{(d)}}K^{(p)}_{p-d}(O_{X_{A},x} \ \mathrm{on} \ x)_{\mathbb{Q}} @>>> \bigoplus\limits_{x \in X ^{(d)}}K^{(p)}_{p-d}(O_{X,x} \ \mathrm{on} \ x)_{\mathbb{Q}} \\
     @VVV @VVV @VVV\\
     0 @. 0 @. 0,
  \end{CD}
\]
}in which the Zariski sheafification of each column is a flasque resolution of $\overline{HC}^{(p-1)}_{p-1}(O_{X_{A}})$,  $K^{(p)}_{p}(O_{X_{A}})_{\mathbb{Q}}$ and $K^{(p)}_{p}(O_{X})_{\mathbb{Q}}$ respectively. The map from the middle column to the left one, denoted $\mathrm{Ch}$, is induced by relative Chern characters from K-theory to cyclic homology, and the map from the middle column to the right one, denoted $\mathrm{Pr}$, is induced by augmentation $A \to k$.

\end{theorem}

Using tensor triangular geometry \cite{Ba1}, Balmer \cite{Ba} defined tensor triangular Chow groups of a tensor triangulated category, which were further explored by Klein \cite{Klein}. By slight modifying Balmer's definition, we proposed Milnor K-theoretic cycles.
\begin{definition}[Definition 3.4 of \cite{Y2}] \label{definition: MilnorKChow}
In notation of Theorem \ref{theorem: Main-Diag}, let $p$ further satisfy that $1 \leq p \leq d$, where $d=\mathrm{dim}(X)$. The $p$-th Milnor K-theoretic cycle groups of $X$ and $X_A$, denoted $Z^{M}_{p}(D^{\mathrm{perf}}(X))$ and $Z^{M}_{p}(D^{\mathrm{perf}}(X_A))$ respectively, are defined to be 
{\footnotesize
\[
  Z^{M}_{p}(D^{\mathrm{perf}}(X)):= \mathrm{Ker}(d_{1,X}^{p,-p}), \ Z^{M}_{p}(D^{\mathrm{perf}}(X_A)):= \mathrm{Ker}(d_{1,X_A}^{p,-p}).
\]
}
The $p$-th Milnor K-theoretic Chow groups of $X$ and $X_A$, denoted by $CH^{M}_{p}(D^{\mathrm{perf}}(X))$ and $CH^{M}_{p}(D^{\mathrm{perf}}(X_A))$ respectively, are defined to be 
{\footnotesize
\[
   CH^{M}_{p}(D^{\mathrm{perf}}(X)):= \dfrac{\mathrm{Ker}(d_{1,X}^{p,-p})}{\mathrm{Im}(d_{1,X}^{p-1,-p})}, \ CH^{M}_{p}(D^{\mathrm{perf}}(X_A)):= \dfrac{\mathrm{Ker}(d_{1,X_A}^{p,-p})}{\mathrm{Im}(d_{1,X_A}^{p-1,-p})}.
\]
}
\end{definition}

The elements of $Z^{M}_{p}(D^{\mathrm{Perf}}(X_{A}))$ are called Milnor K-theoretic cycles. The reason why we use the kernel of $d_{1,X_{A}}^{p,-p}$ to define $Z^{M}_{p}(D^{\mathrm{Perf}}(X_{A}))$ is explained in section 2.2 of \cite{Y3}, where $A$ is the ring of dual numbers $k[t]/(t^{2})$.

To explain that the above definitions are a honest generalization of the classical cycle group $Z^{p}(X)$ and Chow group $CH^{p}(X)$, we recall that
\begin{theorem} [Theorem 3.16 of \cite{Y2}] \label{theorem:CycleGp-Agree}
For $X$ a smooth projective variety over a field $k$ of characteristic zero, there exists the following identifications
{\footnotesize
\[
 Z_{p}^{M}(D^{\mathrm{perf}}(X))= Z^{p}(X)_{\mathbb{Q}}, \ CH^{M}_{p}(D^{\mathrm{perf}}(X)) = CH^{p}(X)_{\mathbb{Q}}.
\]
}

\end{theorem}

In fact, the right column of the diagram in Theorem \ref{theorem: Main-Diag} agrees with the following complex of Milnor K-theory studied by Soul\'e \cite{Soule}
{\footnotesize
\[
0 \to K^{M}_{p}(k(X))_{\mathbb{Q}} \to \bigoplus\limits_{x \in X ^{(1)}}K^{M}_{p-1}(k(x))_{\mathbb{Q}} \to \cdots \to \bigoplus\limits_{x \in X ^{(p)}}K^{M}_{0}(k(x))_{\mathbb{Q}} \to 0.
\]
}This is the key to prove Theorem \ref{theorem:CycleGp-Agree}.

The Milnor K-theoretic Chow groups $CH^{M}_{p}(D^{\mathrm{perf}}(X_{A}))$ agrees with cohomological Chow group $H^{p}(X, K^{M}_{p}(O_{X_A}))_{\mathbb{Q}}$ as follows.
\begin{theorem}  \label{theorem:Agree Milnor Chow}
With notation as above, there are isomorphisms
{\footnotesize
\begin{equation}
 CH^{M}_{p}(D^{\mathrm{perf}}(X_{A})) = H^{p}(X, K^{(p)}_{p}(O_{X_A}))_{\mathbb{Q}}=H^{p}(X, K^{M}_{p}(O_{X_A}))_{\mathbb{Q}}.
\end{equation}
}
\end{theorem}

The first isomorphism follows from Theorem \ref{theorem: Main-Diag} and the second one is from the isomorphism $K^{M}_{p}(O_{X_A})_{\mathbb{Q}}=K^{(p)}_{p}(O_{X_A})_{\mathbb{Q}}$. This answers Question \ref{q:GG's ques}, i.e., extends Bloch formula (1.1) from $X$ to $X_A$. When $A=k[t]/(t^{j})$, the above isomorphisms (2.12) were proved in Theorem 3.17 of \cite{Y2}.

\section{Local Hilbert function and Chow groups}

In this section, we first construct a map from local Hilbert functor to local homology in (3.2). Then, with suitable assumptions, we use this map to answer Question \ref{q:Bloch's ques} in Theorem \ref{t: transf-Hilb-Chow}.

In notation of Setting \ref{s:setting}, let $Y \subset X$ be a closed irreducible subvariety of codimension $p$. 

From now on, we fix the integer $p$. 

\begin{definition} \label{d:Hilbert}
The local Hilbert functor $\mathrm{\mathbb{H}ilb}$ is a functor on the category $Art_{k}$
\begin{align*}
\mathrm{\mathbb{H}ilb}: \ A \longrightarrow \mathrm{\mathbb{H}ilb}(A),
\end{align*}
where $A \in Art_{k}$ and $\mathrm{\mathbb{H}ilb}(A)$ denotes the set of infinitesimal embedded deformations of $Y$ in $X_A$. 
\end{definition}
This functor $\mathrm{\mathbb{H}ilb}$ had been studied intensively in literature, including \cite{Hart2, Sernesi, Vistoli}. We connect it with K-theory in the following.

Let $y$ be the generic point of $Y$. We denote by $K_{0}(O_{X_{A},y} \ \mathrm{on} \ y)$ Grothendieck group of the triangulated category $D^{b}(O_{X_{A},y} \ \mathrm{on} \ y)$, which is the derived category of perfect complexes of $ O_{X_{A},y}$-modules with homology supported on the closed point $y \in \mathrm{Spec}(O_{X,y})$.

The closed subvariety $Y$ is generically given by a regular sequence $\{f_{1}, \cdots, f_{p}\}$ of $O_{X,y}$. For any $Y' \in \mathrm{\mathbb{H}ilb}(A)$, $Y'$ is generically given by a regular sequence $\{f^{A}_{1}, \cdots, f^{A}_{p}\}$ of $O_{X_{A},y}$. Let $L^{A}_{\bullet}$ be the Koszul complex of the regular sequence $\{f^{A}_{1}, \cdots, f^{A}_{p}\}$, we consider the complex $L^{A}_{\bullet}$ as an element of $K_{0}(O_{X_{A},y} \ \mathrm{on} \ y)_{\mathbb{Q}}$.

Adams operations $\psi^{m}$ for K-theory of perfect complexes defined in \cite{GilletSoule} has the following property.
\begin{lemma} [Prop 4.12  of \cite{GilletSoule}]  \label{Lemma: GilletSoule}
Adams operations $\psi^{m}$ on $L^{A}_{\bullet}$ satisfies that
\[
\psi^{m}(L^{A}_{\bullet})  = m^{p}L^{A}_{\bullet}.
\]
\end{lemma}

Let $K^{(p)}_{0}(O_{X_{A},y} \ \mathrm{on} \ y)$ be the eigenspace of $\psi^{m}=m^{p}$. The above Lemma implies that $L^{A}_{\bullet} \in K^{(p)}_{0}(O_{X_{A},y} \ \mathrm{on} \ y)_{\mathbb{Q}}$.

\begin{definition} \label{Definition:Hilb-Map-K}
With notation as above, one defines a set-theoretic map 
\begin{align}
\alpha_{A}: \  \mathrm{\mathbb{H}ilb}(A) &  \longrightarrow K^{(p)}_{0}(O_{X_{A},y} \ \mathrm{on} \ y)_{\mathbb{Q}} \\
Y^{'} \ &  \longrightarrow   L^{A}_{\bullet}.    \notag
\end{align}

\end{definition}

It is interesting to determine whether $\alpha_{A}(Y')$ is a Milnor K-theoretic cycle (in the sense of Definition \ref{definition: MilnorKChow}) or not.
\begin{question} \label{q:K-cycle or not}
With notation as above, is it true that 
\[
\alpha_{A}(Y') \in Z^{M}_{p}(D^{\mathrm{Perf}}(X_{A}))?
\]
In other words, is it true that $d_{1,X_A}^{p,-p} \circ \alpha_{A}(Y')=0$, where $d_{1,X_A}^{p,-p}$ is the differential of the diagram in Theorem \ref{theorem: Main-Diag}?
\end{question}

The subtlety of this question is the appearance of negative K-group $K_{-1}(O_{X_{A},x} \ \mathrm{on} \ x)$ in the diagram in Theorem \ref{theorem: Main-Diag}, which may not vanish. We refer to \cite{CHSW, KST} for recent progress on Weibel's vanishing conjecture of negative K-theory.

In the diagram in Theorem \ref{theorem: Main-Diag}, $K_{-1}(O_{X,x} \ \mathrm{on} \ x)=K_{-1}(k(x))=0$, this implies that $d_{1,X}^{p,-p} \circ \mathrm{Pr} \circ \alpha_{A}(Y')=0$. Since the diagram in Theorem \ref{theorem: Main-Diag} is split, Question  \ref{q:K-cycle or not} is equivalent to the following one.
 \begin{question} \label{q:K-cycle or not-2}
 Let $\mathrm{Ch} \circ \alpha_{A}$ be the composition
 {\footnotesize
  \begin{align}
  \mathrm{\mathbb{H}ilb}(A)  \xrightarrow{\alpha_{A}} K^{(p)}_{0}(O_{X_{A},y} \ \mathrm{on} \ y)_{\mathbb{Q}} \xrightarrow{\mathrm{Ch}} 
  H^{p}_{y}(\overline{HC}^{(p-1)}_{p-1}(O_{X_{A}})),
  \end{align}
}does the image $\mathrm{Ch} \circ \alpha_{A}(Y')$ lie in the kernel of $\partial_{1,X_{A}}^{p,-p}$, where $\mathrm{Ch}$ and $\partial_{1,X_{A}}^{p,-p}$ are maps of the diagram in Theorem \ref{theorem: Main-Diag}?
\end{question}

It is known that $\mathrm{Ch} \circ \alpha_{A}(Y')$ does not always lie in the kernel of $\partial_{1,X_{A}}^{p,-p}$, see Example 4.4 of \cite{Y1}. Hence, $\alpha_{A}(Y')$ is not a Milnor K-theoretic cycle in general.

In the rest of this section, we strength the situation of Setting \ref{s:setting} as follows.
\begin{setting} \label{s:set2}
In notation of Setting \ref{s:setting}, we further assume that $Y \subset X$ is a locally complete intersection. There exists a finite open affine covering $\{U_{i} \}_{i \in I}$ of $X$ such that $Y \cap U_{i}$ is given by a regular sequence $f_{1}, \cdots, f_{p}$ of $O_{X}(U_{i})$.

Let $GArt_{k}\subset Art_{k}$ denote the subcategory of $Art_{k}$ whose objects are also graded $k$-algebras $A=\oplus_{m \geq 0}A_{m}$ such that $A_{0}=k$.
\end{setting}

In this setting, for $A \in GArt_{k}$ and for $Y' \in \mathrm{\mathbb{H}ilb}(A)$, we will prove that $\mathrm{Ch} \circ \alpha_{A}(Y')$ lies in the kernel of $\partial_{1,X_{A}}^{p,-p}$, which yields that the image $\alpha_{A}(Y')$ is a Milnor K-theoretic cycle. 

Let $U_{i}=\mathrm{Spec}(R) \subset X$ be open affine, for $A \in GArt_{k}$, we note that $R \otimes_{k}A= \oplus_{m \geq 0}(R \otimes_{k}A_{m})$ is a graded $k$-algebra with $R \otimes_{k}A_{0}=R \otimes_{k}k=R$. Since $\mathbb{Q} \subset k$, $R \otimes_{k}A$ can be also considered as a graded $\mathbb{Q}$-algebra. By Goodwillie \cite{GoodW}, the SBI sequence (defined over $\mathbb{Q}$) broke into short exact sequence
{\footnotesize
\begin{equation}
0 \to \overline{HC}^{(l-1)}_{l-1}(R \otimes_{k}A) \xrightarrow{\mathrm{B}} \overline{HH}^{(l)}_{l}(R \otimes_{k}A)  \xrightarrow{\mathrm{I}} \overline{HC}^{(l)}_{l}(R \otimes_{k}A) \to 0,
\end{equation}
}where $R \otimes_{k}A$ is considered as a graded $\mathbb{Q}$-algebra and $l$ is any positive integer, $\overline{HC}^{(l-1)}_{l-1}(R \otimes_{k}A)$ is defined to be the kernel of $HC^{(l-1)}_{l-1}(R \otimes_{k}A)\to HC^{(l-1)}_{l-1}(R)$, $\overline{HH}^{(l)}_{l}(R \otimes_{k}A)$ and $\overline{HC}^{(l)}_{l}(R \otimes_{k}A)$ are defined similarly. This sequence is useful in computing cyclic homology and K-theory, for example, see Geller, Reid and Weibel \cite{GRW, GW}. The following short exact sequence
{\footnotesize
\begin{equation}
0 \to \overline{HC}^{(l-1)}_{l-1}(O_{X_{A}}) \xrightarrow{\mathrm{B}} \overline{HH}^{(l)}_{l}(O_{X_{A}})  \xrightarrow{\mathrm{I}} \overline{HC}^{(l)}_{l}(O_{X_{A}}) \to 0,
\end{equation}
}is a sheaf version of (3.3).

For each integer $q$ satisfying that $1 \leq q \leq d$, where $d=\mathrm{dim}(X)$, and for $x \in X^{(q)}$, there is a long exact sequence associated to (3.4)
{\footnotesize
\begin{align}
& \cdots \to H^{*}_{x}(\overline{HC}^{(l-1)}_{l-1}(O_{X_{A}})) \to H^{*}_{x}(\overline{HH}^{(l)}_{l}(O_{X_{A}})) \to H^{*}_{x}(\overline{HC}^{(l)}_{l}(O_{X_{A}}))  \\
& \to H^{*+1}_{x}(\overline{HC}^{(l-1)}_{l-1}(O_{X_{A}}))  \to \cdots. \notag
\end{align}
}

By Corollary \ref{c: HC-CM} and Lemma \ref{l: HH-CM}, the sheaves of sequence (3.4) are Cohen-Macaulay, so the sequence (3.5) is indeed a short exact sequence
{\footnotesize
\begin{align}
 0 \to  H^{q}_{x}(\overline{HC}^{(l-1)}_{l-1}(O_{X_{A}})) \xrightarrow{\mathrm{B}} H^{q}_{x}(\overline{HH}^{(l)}_{l}(O_{X_{A}})) \xrightarrow{\mathrm{I}} H^{q}_{x}(\overline{HC}^{(l)}_{l}(O_{X_{A}})) \to 0.   
\end{align}
}

We recall that $y$ is the generic point of $Y$, $y\in X^{(p)}$. To investigate Question \ref{q:K-cycle or not-2}, we want to describe the composition
\begin{align}
&\mathrm{\mathbb{H}ilb}(A)  \xrightarrow{\alpha_{A}} K^{(p)}_{0}(O_{X_{A},y} \ \mathrm{on} \ y)_{\mathbb{Q}} \\ & \xrightarrow{\mathrm{Ch}} H^{p}_{y}(\overline{HC}^{(p-1)}_{p-1}(O_{X_{A}})) \xrightarrow{\mathrm{B}} H^{p}_{y}( \overline{HH}^{(p)}_{p}(O_{X_{A}})), \notag
\end{align}
where $\mathrm{B}$ is the injective map in (3.6) (let $l=p$, $q=p$ and $x=y$).

Let $\overline{\Omega}_{X_{A} / \mathbb{Q}}^{p}$ be the kernel of $\Omega_{X_{A} / \mathbb{Q}}^{p} \to \Omega_{X/ \mathbb{Q}}^{p}$, by Lemma \ref{lemma: l-l-omega}, there is an isomorphism
\begin{equation}
\overline{HH}^{(p)}_{p}(O_{X_{A}}) = \overline{\Omega}_{X_{A} / \mathbb{Q}}^{p}.
\end{equation}
Then we write the above composition (3.7) as
\begin{align}
&\mathrm{\mathbb{H}ilb}(A)  \xrightarrow{\alpha_{A}} K^{(p)}_{0}(O_{X_{A},y} \ \mathrm{on} \ y)_{\mathbb{Q}} \\ &\xrightarrow{\mathrm{Ch}}  H^{p}_{y}(\overline{HC}^{(p-1)}_{p-1}(O_{X_{A}})) \xrightarrow{\mathrm{B}} H^{p}_{y}( \overline{\Omega}_{X_{A} / \mathbb{Q}}^{p}), \notag
\end{align}
and describe it in the following.

We first use a construction of Ang\'eniol and Lejeune-Jalabert \cite{ALJ} to describe the composition $\mathrm{B} \circ \mathrm{Ch}$\footnote{Analogous descriptions were given in \cite{Y1} (Section 3) and \cite{Y4} (Section 2) by using Ang\'eniol and Lejeune-Jalabert's method, where $A$ is a truncated polynomial $k[t]/(t^{j})$.}. An element of $K^{(p)}_{0}(O_{X_{A},y} \ \mathrm{on} \ y)_{\mathbb{Q}}$ is represented by a strict perfect complex $L_{\bullet}$
{\small
\[
 \begin{CD}
  0 @>>> L_{n} @>M_{n}>> L_{n-1} @>M_{n-1}>>  \dots @>M_{2}>> L_{1} @>M_{1}>> L_{0} @>>> 0,
 \end{CD}
\]
}where each $L_{j}$ is a free $O_{X_A,y}$-modules of finite rank, each $M_{j}$ is a matrix with entries in $O_{X_A,y}$ and the homology of $L_{\bullet}$ is supported on $y$.

\begin{definition} [page 24 in \cite{ALJ}] \label{definition: local cycles}
The local fundamental class attached to this perfect complex $L_{\bullet}$ is defined to be the following collection
\[
 [L_{\bullet}]_{loc}=\{\frac{1}{p!}dM_{j}\circ dM_{j+1}\circ \dots \circ dM_{j+p-1}\}, j =  1, 2, \cdots,
\]
where $d=d_{\mathbb{Q}}$ and each $dM_{j}$ is the matrix of differentials. In other words,
\[
 dM_{j} \in \mathrm{Hom}(L_{j},L_{j-1}\otimes \Omega_{O_{X_A,y} / \mathbb{Q}}^{1}).
\]
\end{definition}

By Lemma 3.1.1 (on page 24) and Definition 3.4 (on page 29) in \cite{ALJ}, the local fundamental class $[L_{\bullet}]_{loc}$ defines a cycle of the complex $\mathcal{H}om( L_{\bullet}, \Omega^{p}_{O_{X_A,y} / \mathbb{Q}}\otimes L_{\bullet})$ and its image (still denoted $[L_{\bullet}]_{loc}$) in $\mathcal{E}XT^{p}(L_{\bullet}, \Omega^{p}_{O_{X_A,y}/ \mathbb{Q}}\otimes L_{\bullet})$, which is the $p$-th cohomology of the complex $\mathcal{H}om(L_{\bullet}, \Omega^{p}_{O_{X_A,y}/ \mathbb{Q}}\otimes L_{\bullet})$, does not depend on the choice of the basis of $L_{\bullet}$.

Since $L_{\bullet}$ is supported on $y$, by the discussion after Definition 2.3.1 on page 98-99 in \cite{ALJ}, there exists a trace map
\[
 \mathrm{Tr}:  \mathcal{E}XT^{p}(L_{\bullet}, \Omega^{p}_{O_{X_A,y}/ \mathbb{Q}}\otimes L_{\bullet}) \to H^{p}_{y}(\Omega^{p}_{X_A/ \mathbb{Q}}).
\]

\begin{definition}   [Definition 2.3.2 on page 99 in \cite{ALJ}] \label{Definition: Newton class}
The image of $[L_{\bullet}]_{loc}$ under the above trace map $\mathrm{Tr}$, denoted $\mathcal{V}^{p}_{L_{\bullet}}$, is called Newton class.
\end{definition}

Grothendieck group of a triangulated category is the monoid of isomorphism objects modulo the submonoid formed from distinguished triangles.
\begin{lemma} [Proposition 4.3.1 on page 113 in \cite{ALJ}] \label{Lemma: Newton class-Well}
The Newton class $\mathcal{V}^{p}_{L_{\bullet}}$ is well-defined on the Grothendieck group $K_{0}(O_{X_A,y} \ \mathrm{on} \ y)$.
\end{lemma}

The morphism $\Omega^{p}_{X_A / \mathbb{Q}} \to \overline{\Omega}^{p}_{X_A / \mathbb{Q}}$ induces a map $
 \varphi: H^{p}_{y}(\Omega^{p}_{X_A / \mathbb{Q}}) \to H^{p}_{y}(\overline{\Omega}^{p}_{X_A / \mathbb{Q}})$.

\begin{definition} \label{Definition: Chern-Newton}
One uses Newton class $\mathcal{V}^{p}_{L_{\bullet}}$ to defines a morphism
\begin{align*}
\rho: K^{(p)}_{0}&(O_{X_{A},y} \ \mathrm{on} \ y)_{\mathbb{Q}} \to  H^{p}_{y}(\Omega^{p}_{X_A / \mathbb{Q}})  \to H^{p}_{y}(\overline{\Omega}^{p}_{X_A / \mathbb{Q}}) \\
& \ \ \ L_{\bullet} \ \ \ \ \ \ \ \ \ \  \longrightarrow  \ \ \ \ \  \mathcal{V}^{p}_{L_{\bullet}} \ \ \  \longrightarrow  \ \ \ \  \varphi(\mathcal{V}^{p}_{L_{\bullet}}).
\end{align*}

\end{definition}

The composition $\mathrm{B} \circ \mathrm{Ch}$ in (3.9) can be described by $\rho$, so $\mathrm{B} \circ \mathrm{Ch} \circ \alpha_{A}$ in (3.9) is given by $\rho \circ \alpha_{A}$. Concretely, in notation of Setting \ref{s:set2}, for any $Y' \in \mathrm{\mathbb{H}ilb}(A)$, $Y'$ is still a locally complete intersection.  In fact, $Y' \cap U_{i}$ is given by a regular sequence $\{f^{A}_{1}, \cdots, f^{A}_{p}\}$ of $O_{X_A}(U_{i})$. 

By considering each $f_{i}$ and $f^{A}_{i}$ as elements of $O_{X,y}$ and $O_{X_{A},y}$ respectively, one has that $Y$ and $Y'$ are generically given by regular sequences $\{f_{1}, \cdots, f_{p}\}$ and $\{f^{A}_{1}, \cdots, f^{A}_{p}\}$ respectively.

Let $F^{A}_{\bullet}$ be the Koszul resolution of $O_{X_{A},y}/(f^{A}_{1}, \cdots, f^{A}_{p})$, which has the form
\[
 0 \to F^{A}_{p} \to \cdots \to F^{A}_{0} \to 0,
\]
where each $F^{A}_{i}=\bigwedge^{i}(O_{X_{A},y})^{\oplus p}$. 

By Definition \ref{Definition:Hilb-Map-K}, $\alpha_{A}(Y') = F^{A}_{\bullet} \in K^{(p)}_{0}(O_{X_{A},y} \ \mathrm{on} \ y)_{\mathbb{Q}}$. The image $\mathrm{B} \circ \mathrm{Ch} \circ \alpha_{A}(Y')$ can be described via Newton class. Concretely, the following diagram
{\tiny
\[
\begin{cases}
 \begin{CD}
  F^{A}_{\bullet} @>>> O_{X_{A},y}/(f^{A}_{1}, \cdots, f^{A}_{p})  \\
  F^{A}_{p}(\cong O_{X_{A},y}) @> [F^{A}_{\bullet}]_{loc}>> F^{A}_{0} \otimes \Omega_{O_{X_{A},y}/ \mathbb{Q}}^{p} (\cong \Omega_{O_{X_{A},y}/ \mathbb{Q}}^{p}), 
 \end{CD}
\end{cases}
\]
}where $[F^{A}_{\bullet}]_{loc}=df^{A}_{1}  \wedge \cdots \wedge df^{A}_{p}$ is the local fundamental class attached to $F^{A}_{\bullet}$, gives an element $\beta^{A}$ in {\small $Ext^{p}(O_{X_{A},y}/(f^{A}_{1}, \cdots, f^{A}_{p}), \Omega_{O_{X_{A},y}/ \mathbb{Q}}^{p})$}. There is an isomorphism
\[
H_{y}^{p}(\Omega_{X_{A}/ \mathbb{Q}}^{p}) = \varinjlim_{n \to \infty}Ext^{p}(O_{X_{A},y}/(f^{A}_{1}, \cdots, f^{A}_{p})^{n}, \Omega_{O_{X_{A},y}/ \mathbb{Q}}^{p}),
\]
the image $[\beta^{A}]$ of $\beta^{A}$ under the limit is the Newton class $\mathcal{V}^{p}_{F^{A}_{\bullet}} \in H_{y}^{p}(\Omega_{X_{A}/ \mathbb{Q}}^{p})$.

Let $F_{\bullet}(f_{1}, \cdots, f_{p})$ be the Koszul resolution of $O_{X,y}/(f_{1}, \cdots, f_{p})$, which has the form
\[
 0 \to F_{p} \to \cdots \to F_{0} \to 0,
\]
where each $F_{i}$ is defined as usually.

For $[F^{A}_{\bullet}]_{loc}=df^{A}_{1}  \wedge \cdots \wedge df^{A}_{p} \in \Omega^{p}_{O_{X_A,y}/ \mathbb{Q}}$, we denote by $\overline{[F^{A}_{\bullet}]}_{loc}$ the image of $[F^{A}_{\bullet}]_{loc}$ under the morphism $\Omega^{p}_{O_{X_A,y}/ \mathbb{Q}} \to \overline{\Omega}^{p}_{O_{X_A,y}/ \mathbb{Q}}$, where $\overline{\Omega}^{p}_{O_{X_A,y}/ \mathbb{Q}}$ is the kernel of $\Omega^{p}_{O_{X_A,y}/ \mathbb{Q}} \to \Omega^{p}_{O_{X,y}/ \mathbb{Q}}$. Concretely, $\overline{[F^{A}_{\bullet}]}_{loc}=df^{A}_{1}  \wedge \cdots \wedge df^{A}_{p}-df_{1}  \wedge \cdots \wedge df_{p}$. The following diagram (denoted $\overline{\beta^{A}}$)
\begin{equation}
\begin{cases}
 \begin{CD}
    F_{\bullet}(f_{1}, \cdots, f_{p})  @>>>   O_{X,y}/(f_{1}, \cdots, f_{p}) \\
  F_{p}(\cong O_{X,y}) @>\overline{[F^{A}_{\bullet}]}_{loc}>> F_{0} \otimes \overline{\Omega}^{p}_{O_{X_A,y}/ \mathbb{Q}}(\cong \overline{\Omega}^{p}_{O_{X_A,y}/ \mathbb{Q}}),
 \end{CD}
 \end{cases}
\end{equation}
defines an element in $Ext^{p}(O_{X,y}/(f_{1}, \cdots, f_{p}),\overline{\Omega}^{p}_{O_{X_A,y}/ \mathbb{Q}})$. There is an isomorphism
\[
H_{y}^{p}(\overline{\Omega}_{X_{A}/ \mathbb{Q}}^{p}) = \varinjlim_{n \to \infty}Ext^{p}(O_{X,y}/(f_{1}, \cdots, f_{p})^{n}, \overline{\Omega}_{O_{X_{A},y}/ \mathbb{Q}}^{p}),
\]
the image $[\overline{\beta^{A}}] \in H_{y}^{p}(\overline{\Omega}^{p}_{X_{A}/ \mathbb{Q}})$ of $\overline{\beta^{A}}$ under the limit is $\mathrm{B} \circ \mathrm{Ch} \circ \alpha_{A}(Y')$.

To summarize, one has
\begin{lemma} \label{l: describe B-Ch-Alpha}
In Setting \ref{s:set2}, for $Y' \in \mathrm{\mathbb{H}ilb}(A)$, the image of $Y'$ under the composition $\mathrm{B} \circ \mathrm{Ch} \circ \alpha_{A}$ in (3.9) can be described by $[\overline{\beta^{A}}]$
\[
\mathrm{B} \circ \mathrm{Ch} \circ \alpha_{A}(Y')=[\overline{\beta^{A}}].
\]
\end{lemma}

Let $l=p$ in the sequence (3.4), the natural map $\mathrm{B}: \overline{HC}^{(p-1)}_{p-1}(O_{X_{A}}) \to \overline{HH}^{(p)}_{p}(O_{X_{A}})$ induces a commutative diagram
{\footnotesize
\begin{equation}
  \begin{CD}
     0 @. 0 \\
      @VVV @VVV \\
      \overline{HC}^{(p-1)}_{p-1}(O_{X_{A},\eta}) @>\mathrm{B}>>  \overline{\Omega}^{p}_{O_{X_{A},\eta}/\mathbb{Q}} \\
     @VVV @VVV \\
      \bigoplus\limits_{x \in X ^{(1)}} H^{1}_{x}(\overline{HC}^{(p-1)}_{p-1}(O_{X_{A}}))  @>\mathrm{B}>> \bigoplus\limits_{x \in X ^{(1)}} H^{1}_{x}(\overline{\Omega}^{p}_{X_{A}/\mathbb{Q}})  \\
    @VVV  @VVV \\
     \dots @>>>  \dots \\ 
      @VVV @VVV \\
     \bigoplus\limits_{x \in X ^{(p)}} H^{p}_{x}(\overline{HC}^{(p-1)}_{p-1}(O_{X_{A}}))  @>\mathrm{B}>> \bigoplus\limits_{x \in X ^{(p)}} H^{p}_{x}(\overline{\Omega}^{p}_{X_{A}/\mathbb{Q}})  \\
     @V \partial_{1,X_{A}}^{p,-p}VV  @V \tilde{\partial}_{1,X_{A}}^{p,-p}VV \\
    \bigoplus\limits_{x \in X ^{(p+1)}} H^{p+1}_{x}(\overline{HC}^{(p-1)}_{p-1}(O_{X_{A}}))  @>\mathrm{B}>> \bigoplus\limits_{x \in X ^{(p+1)}} H^{p+1}_{x}(\overline{\Omega}^{p}_{X_{A}/\mathbb{Q}})  \\
      @VVV  @VVV\\
     \dots @>>> \dots \\ 
     @VVV @VVV \\
      \bigoplus\limits_{x \in X ^{(d)}} H^{d}_{x}(\overline{HC}^{(p-1)}_{p-1}(O_{X_{A}}))  @>\mathrm{B}>> \bigoplus\limits_{x \in X ^{(d)}} H^{d}_{x}(\overline{\Omega}^{p}_{X_{A}/\mathbb{Q}})  \\
     @VVV @VVV \\
     0 @. 0,
  \end{CD}
\end{equation}
}where the two columns are Cousin complexes of $\overline{HC}^{(p-1)}_{p-1}(O_{X_{A}})$ and $\overline{HH}^{(p)}_{p}(O_{X_{A}})$ respectively and we use (3.8) to identify $\overline{HH}^{(p)}_{p}(O_{X_{A}})$ with $\overline{\Omega}^{p}_{X_{A}/\mathbb{Q}}$.

\begin{lemma} \label{l:kernel-HH}
With notation as above, for $[\overline{\beta^{A}}] \in H_{y}^{p}(\overline{\Omega}^{p}_{X_{A}/ \mathbb{Q}})$, where $\overline{\beta^{A}}$ is (3.10), one has
\[
\tilde{\partial}_{1,X_{A}}^{p,-p} ([\overline{\beta^{A}}])=0,
\]
where $\tilde{\partial}_{1,X_{A}}^{p,-p}$ is the differential of the right column of diagram (3.11). In other words, for $Y' \in \mathrm{\mathbb{H}ilb}(A)$, the image $\mathrm{B} \circ \mathrm{Ch} \circ \alpha(Y')$ in (3.9) lies in the kernel of $\tilde{\partial}_{1,X_{A}}^{p,-p}$
\[
\tilde{\partial}_{1,X_{A}}^{p,-p} \circ \mathrm{B} \circ \mathrm{Ch} \circ \alpha(Y')=0.
\]
\end{lemma}

\begin{proof}
In notation of Setting \ref{s:set2}, by shrinking $U_{i}$, we assume that $O_{X}(U_{i})$ is local. The regular sequence $\{f_{1}, \cdots, f_{p}\}$ can be extended to a system of parameter $\{f_{1}, \cdots, f_{p}, f_{p+1}, \cdots, f_{d}\}$ of the regular local ring $O_{X}(U_{i})$. The prime ideals $Q_{j}:=(f_{1}, \cdots, f_{p}, f_{j})$, where $j= p+1, \cdots, d$, define generic points $z_{j} \in X^{(p+1)}$. In the following, to check $\tilde{\partial}_{1,X_{A}}^{p,-p} \circ \mathrm{B} \circ \mathrm{Ch} \circ \alpha(Y')=0$, we consider the prime $Q_{p+1}=(f_{1}, \cdots, f_{p}, f_{p+1})$ which defines the generic point $z_{p+1}$, other cases work similarly.  

Let $Q=(f_{1}, \cdots, f_{p})$ be the prime ideal which defines the generic point (of $Y$) $y \in X^{(p)}$, then $O_{X,y}=(O_{X,z_{p+1}})_{Q}$. Then $\overline{\beta^{A}}$ (cf. (3.10)) can be rewritten as
{\footnotesize
\begin{equation*}
\begin{cases}
 \begin{CD}
   F_{\bullet}(f_{1}, f_{2},\cdots, f_{p}) @>>> (O_{X,z_{p+1}})_{Q}/(f_{1}, f_{2}, \cdots,  f_{p}) \\
  F_{p}(\cong (O_{X,z_{p+1}})_{Q}) @> \dfrac{f_{p+1}}{f_{p+1}}\overline{[F^{A}_{\bullet}]}_{loc}>> F_{0} \otimes \overline{\Omega}_{(O_{X_{A},z_{p+1}})_{Q} / \mathbb{Q}}^{p}(\cong \overline{\Omega}_{(O_{X_{A},z_{p+1}})_{Q}/ \mathbb{Q}}^{p}).
 \end{CD}
 \end{cases}
\end{equation*}
}Here $\overline{\Omega}_{(O_{X_{A},z_{p+1}})_{Q}/ \mathbb{Q}}^{p}$ is the kernel of $\Omega_{(O_{X_{A},z_{p+1}})_{Q}/ \mathbb{Q}}^{p} \to \Omega_{(O_{X,z_{p+1}})_{Q}/ \mathbb{Q}}^{p}$, and $F_{\bullet}(f_{1}, f_{2},\cdots, f_{p})$ is of the form
{\footnotesize
\[
 \begin{CD}
  0 @>>> F_{p} @>>> F_{p-1} @>>>  \dots @>>> F_{1} @>>> F_{0},
 \end{CD}
\]
}where each $F_{i}=\bigwedge^{i}((O_{X,z_{p+1}})_{Q})^{\oplus p}$. Since $f_{p+1} \notin Q=(f_{1}, \cdots, f_{p})$, $f_{p+1}^{-1}$ exists in $(O_{X,z_{p+1}})_{Q}$, we can write $\overline{[F^{A}_{\bullet}]}_{loc}=\dfrac{f_{p+1}}{f_{p+1}}\overline{[F^{A}_{\bullet}]}_{loc}$.

The image $\tilde{\partial}_{1,X_{A}}^{p,-p}(\overline{\beta}^{A})$ is represented by the following diagram (denoted $\gamma$)
{\footnotesize
\[
\begin{cases}
 \begin{CD}
   F_{\bullet}(f_{1}, f_{2},\cdots, f_{p}, f_{p+1}) @>>> O_{X, z_{p+1}}/(f_{1}, f_{2}, \cdots,  f_{p}, f_{p+1}) \\
  F_{p+1}(\cong O_{X,z_{p+1}}) @>f_{p+1}\overline{[F^{A}_{\bullet}]}_{loc}>> F_{0} \otimes \overline{\Omega}_{O_{X_{A},z_{p+1}}/ \mathbb{Q}}^{p}(\cong \overline{\Omega}_{O_{X_{A},z_{p+1}}/ \mathbb{Q}}^{p}),
 \end{CD}
 \end{cases}
\]
}where $\overline{\Omega}_{O_{X_{A},z_{p+1}}/ \mathbb{Q}}^{p}$ is the kernel of $\Omega_{O_{X_{A},z_{p+1}}/ \mathbb{Q}}^{p} \to \Omega_{O_{X,z_{p+1}}/ \mathbb{Q}}^{p}$ and the complex $F_{\bullet}(f_{1}, f_{2},\cdots, f_{p}, f_{p+1})$ is of the form
{\footnotesize
\[
 \begin{CD}
  0 @>>>  \bigwedge^{p+1}(O_{X, z_{p+1}})^{\oplus p+1} @>M_{p+1}>> \bigwedge^{p}(O_{X, z_{p+1}})^{\oplus p+1} @>>>  \cdots.
 \end{CD}
\]
}

Let $\{e_{1}, \cdots, e_{p+1} \}$ be a basis of $(O_{X, z_{p+1}})^{\oplus p+1} $, the map $M_{p+1}$ is 
{\footnotesize
\[
 e_{1}\wedge \cdots \wedge e_{p+1}  \to \sum^{p+1}_{j=1}(-1)^{j}f_{j}e_{1}\wedge \cdots \wedge \hat{e_{j}} \wedge \cdots e_{p+1},
\]
}where $\hat{e_{j}}$ means to omit $e_{j}$. Since $f_{p+1}$ appears in $M_{p+1}$, one has that
{\footnotesize
\[
\gamma=0 \in Ext_{O_{X,z_{p+1}}}^{p+1}(O_{X, z_{p+1}}/(f_{1}, f_{2}, \cdots,  f_{p}, f_{p+1}), \overline{\Omega}_{O_{X,z_{p+1}}/ \mathbb{Q}}^{p}).
\]
}Hence, $\tilde{\partial}_{1,X_{A}}^{p,-p}(\overline{\beta^{A}})=0$.

\end{proof}

The commutativity of diagram (3.11) yields that $\mathrm{B} \circ \partial_{1,X_{A}}^{p,-p} \circ \mathrm{Ch} \circ \alpha(Y')=\tilde{\partial}_{1,X_{A}}^{p,-p} \circ \mathrm{B} \circ \mathrm{Ch} \circ \alpha(Y')=0$. Each $\mathrm{B}$ map in diagram (3.11) is injective (see the exact sequence (3.6)), so $\partial_{1,X_{A}}^{p,-p} \circ \mathrm{Ch} \circ \alpha(Y')=0$. This answers Question \ref{q:K-cycle or not-2} in Setting \ref{s:set2}. Equivalently, it answers Question \ref{q:K-cycle or not} in Setting \ref{s:set2}.
\begin{theorem} \label{t:AlphaCycle}
In Setting \ref{s:set2}, for any $A \in GArt_{k}$ and for $Y' \in \mathrm{\mathbb{H}ilb(A)}$, $\alpha_{A}(Y')$ is a Milnor K-theoretic cycle. 
\end{theorem}

The Milnor K-theoretic cycle $\alpha_{A}(Y')$ defines an element of Milnor K-theoretic Chow group (defined in Definition \ref{definition: MilnorKChow}), which further gives an element of the cohomological Chow group $CH^{p}(X,K^{M}_{p}(O_{X_A}))_{\mathbb{Q}}$ by Theorem \ref{theorem:Agree Milnor Chow}, denoted $[\alpha_{A}(Y')]$.
There is a set-theoretic map
\begin{align}
 \mathrm{\mathbb{H}ilb}(A) & \to  \widetilde{CH}^{p}(A) \\  
  \ Y' & \to [\alpha_{A}(Y')],   \notag
\end{align}
where $\widetilde{CH}^{p}(A)=CH^{p}(X,K^{M}_{p}(O_{X_A}))_{\mathbb{Q}}$, see (1.2) on page 1. 

Let $f: C \to A$ be a morphism in the category $GArt_{k}$, there exists a commutative diagram of sets
\begin{equation*}
\begin{CD}
\mathrm{\mathbb{H}ilb}(C) @>\alpha_{C}\mathrm{(3.1)}>>   K^{(p)}_{0}(O_{X_{C},y} \ \mathrm{on} \ y)_{\mathbb{Q}} \\
@Vf_{H}VV @Vf_{K}VV \\
\mathrm{\mathbb{H}ilb}(A) @>\alpha_{A}\mathrm{(3.1)}>> K^{(p)}_{0}(O_{X_{A},y} \ \mathrm{on} \ y)_{\mathbb{Q}},  \\
\end{CD}
\end{equation*}
where $f_{H}$ and $f_{K}$ are induced by $f$ respectively. Since $Y$ is a locally complete intersection, this square can be straightforwardly checked. This induces a commutative diagram of sets
\begin{equation*}
\begin{CD}
\mathrm{\mathbb{H}ilb}(C) @>(3.12)>>   \widetilde{CH}^{p}(C) \\
@VVV @VVV \\
\mathrm{\mathbb{H}ilb}(A) @>(3.12)>> \widetilde{CH}^{p}(A).  \\
\end{CD}
\end{equation*}
We deduce that
\begin{theorem} \label{t: transf-Hilb-Chow}
In Setting \ref{s:set2}, there exists a natural transformation between functors on $GArt_{k}$
\begin{equation*}
\mathbb{T}:  \mathrm{\mathbb{H}ilb} \to \widetilde{CH}^{p},
\end{equation*}
which is defined to be, for any $ A \in GArt_{k}$, $\mathbb{T}(A)$ is (3.12).

\end{theorem}

This answers Bloch's Question \ref{q:Bloch's ques} in Setting \ref{s:set2}. \\

\textbf{Acknowledgements.} This paper is a follow-up to joint work with Benjamin Dribus and Jerome William Hoffman \cite{DHY}, the author thanks both of them and Marco Schlichting for many discussions. He thanks Spencer Bloch \cite{Bl5} for sharing ideas and also thanks Bangming Deng, Phillip Griffiths, Luc Illusie, Marc Levine, Kefeng Liu, Christophe Soul\'e and Jan Stienstra for their suggestions and/or comments on related work \cite{Y1, Y2, Y3, Y4}.


\begin{thebibliography}{99}

 \bibitem{ALJ}
 B. Ang\'eniol and M. Lejeune-Jalabert, \emph{Calcul diff\'erentiel et classes caract\'eristiques en g\'eom\'etrie alg\'ebrique}, (French) [Differential calculus and characteristic classes in algebraic geometry] With an English summary, Travaux en Cours [Works in Progress], 38, Hermann, Paris, 1989.



 
  \bibitem{Ba1}
 P. Balmer, \emph{Tensor triangular geometry}, Proceedings of the International Congress of Mathematicians, Volume II, 85-112, Hindustan Book Agency, New Delhi, 2010.
 
 \bibitem{Ba}
 P. Balmer,  \emph{Tensor triangular Chow groups}, Journal of Geometry and Physics 72 (2013), 3-6. 
 
 
\bibitem{Bl1-tan}
S. Bloch, \emph{On the Tangent Space to Quillen $K$-theory}, Lecture Notes in Mathematics 341 (1972), 205-210. 

 
 \bibitem{Bl2-Annals}
S. Bloch, \emph{$K_2$ and Algebraic Cycles}, Annals of Mathematics 99 (2) (1974), 349-379.

\bibitem{Bl3}
 S. Bloch, \emph{$K_{2}$ of Artinian $\mathbb{Q}$-algebras, with application to algebraic cycles}, Comm. Algebra 3 (1975), 405-428.


 
 \bibitem{Bl4}
 S. Bloch, \emph{Lectures on algebraic cycles}, Second edition, New Mathematical Monographs, 16, Cambridge University Press, Cambridge, 2010, xxiv+130 pp. ISBN: 978-0-521-11842-2.
 
 \bibitem{Bl5}
  S. Bloch, \emph{Private discussions at Tsinghua University},  Beijing, Spring semester (2015) and Fall semester (2016). 
 
 \bibitem{BEK1}
 S. Bloch, H. Esnault and M. Kerz,  \emph{p-adic deformation of algebraic cycle classes}, Invent. Math. 195 (3) (2014), 673-722. 

 \bibitem{BEK2}
 S. Bloch, H. Esnault and M. Kerz,  \emph{Deformation of algebraic cycle classes in characteristic zero}, Algebraic Geometry 1 (3) (2014), 290-310.



\bibitem{BKL}
S. Bloch, A. Kas, D. Lieberman, \emph{Zero cycles on surfaces with $p_{g}=0$}, Compositio Math. 33 (2) (1976), 135-145.


\bibitem{BlochOgus}
S. Bloch and A. Ogus, \emph {Gersten's conjecture and the homology of schemes}, Ann. Sci. {\'E}cole Norm. Sup. 7 (4) (1974), 181-201. 
 
 \bibitem{BS}
S. Bloch, V. Srinivas, \emph{Remarks on correspondences and algebraic cycles}, Amer. J. of Math. 105 (1983), 1235-1253.
 
  \bibitem{Cath}
J.-L. Cathelineau, \emph{$\lambda$-Structures in Algebraic $K$-Theory and Cyclic Homology}, $K$-Theory 4 (1991), 591--606. 


\bibitem{CTHK}
J.-L. Colliot-Th\'{e}l\`{e}ne, R. T. Hoobler and B. Kahn, \emph{The Bloch-Ogus-Gabber Theorem}, Algebraic K-theory (Toronto, ON, 1996), 31Ð-94, Fields Inst. Commun., 16, Amer. Math. Soc., Providence, RI, 1997. 
  
 
 
\bibitem{CHW}
G. Corti\~nas, C. Haesemeyer and C. Weibel, \emph{Infinitesimal cohomology and the Chern character
to negative cyclic homology}, Math. Annalen 344 (2009), 891-922.

\bibitem{CHSW}
G. Corti\~nas, C. Haesemeyer, M. Schlichting and C. Weibel, \emph{Cyclic homology, cdh-cohomology and negative K-theory}, Ann. of Math. 167 (2) (2008), 549-573.



\bibitem{DHY}
B. Dribus, J. W. Hoffman and S. Yang,   \emph{Tangents to Chow Groups: on a question of Green-Griffiths},  Bollettino dell'Unione Matematica Italiana 11 (2018), 205-244.


\bibitem{Gabber}
O. Gabber, \emph{Gersten's conjecture for some complexes of vanishing cycles}, Manuscr. Math. 85 (1994), 323-343.


\bibitem{GRW}
S. Geller, L. Reid and C. Weibel, \emph{The cyclic homology and K-theory of curves}, J. reine angew. Math. 393 (1989), 39-90.



\bibitem{GW}
S. Geller and C. Weibel, \emph{Hodge Decompositions of Loday symbols in K-theory and cyclic homology}, K-theory 8 (1994), 587--632.


\bibitem{GilletSoule}
H. Gillet and C. Soul\'e, \emph{Intersection theory using Adams operations}, Invent. Math. 90 (2) (1987), 243-277.


\bibitem{GoodW}  
 T. G. Goodwillie, \emph{Cyclic homology, derivations and the free loop space}, Topology 24 (1985), 187-215.


  \bibitem{Good}  
 T. G. Goodwillie, \emph{Relative algebraic K-theory and cyclic homology}, Ann. of Math. 124 (2) (1986), 347-402.    
 
 

\bibitem{Gray}
D. Grayson, \emph{Universal exactness in algebraic K-theory}, J. Pure Appl. Algebra 36 (1985), 139--141.

 
 
 



\bibitem{GGChow}
 M. Green and P. Griffiths, \emph{Formal deformation of Chow groups}, The legacy of Niels Henrik Abel. (2004), 467-509, Springer, Berlin.


\bibitem{GGtangentspace}
M. Green and P. Griffiths, \emph{On the Tangent space to the space of algebraic cycles on a smooth algebraic variety}, Annals of Math Studies, 157, Princeton University Press, Princeton, NJ, 2005, vi+200 pp. ISBN: 0-681-12044-7.


\bibitem{Hart1}
R. Hartshorne,  \emph{Residues and duality}, Lecture Notes in Mathematics, 20, Springer-Verlag, Berlin-New York 1966 vii+423 pp.

\bibitem{Hart2}
R. Hartshorne, \emph{Deformation theory}, Graduate Texts in Mathematics, 257, Springer, New York, 2010, viii+234 pp. ISBN: 978-1-4419-1595-5. 
 
 \bibitem{Hu}
W. Hu, \emph{The Chow group of zero cycles for the quotient of certain Calabi-Yau varieties}, Houston J. Math. 38 (1) (2012), 55-67.
    
\bibitem{Keller}
B. Keller, \emph{On the cyclic homology of ringed spaces and schemes}, Doc. Math.  3 (1998), 231-259.

\bibitem{Keller2}
B. Keller, \emph{On the cyclic homology of exact categories}, J. Pure Appl. Algebra 136 (1) (1999), 1-56.

\bibitem{Kerz}
M. Kerz, \emph{The Gersten conjecture for Milnor K-theory}, Invent. math. 175 (2009), 1-33.

\bibitem{KST}
M. Kerz, F. Strunk and G. Tamme, \emph{Algebraic K-theory and descent for blow-ups}, Invent. Math. 211 (2) (2018), 523-577.




\bibitem{Klein}
S. Klein, \emph{Chow groups of tensor-triangulated categories}, J. Pure Appl. Algebra 220 (4) (2016), 1343-1381.
 
 \bibitem{Levine}
M. Levine,
\emph{Lambda-operations, $K$-theory and motivic cohomology}, Algebraic K-theory (Toronto, ON, 1996), 131-184, Fields Inst. Commun., 16, Amer. Math. Soc., Providence, RI, 1997.

\bibitem{Lewis1}
J. Lewis, \emph{Towards a generalization of Mumford's theorem}, J. Math. Kyoto Univ. 29 (1989), 267-272.

\bibitem{Lewis2}

J. Lewis, \emph{A generalization of Mumford's theorem II}, Illinois Journal of Mathematics. 39 (1995), 288-304.

\bibitem{Loday}
 J.-L. Loday, \emph{Cyclic homology}, Grundlehren der Mathematischen Wissenschaften
[Fundamental Principles of Mathematical Sciences], 301, Springer-Verlag, Berlin,
1992. Appendix E by Mar\'ia O. Ronco.





 \bibitem{MaazenStien}
H. Maazen and J. Stienstra, 
\emph {A presentation of $K_2$ of split radical pairs}, J. Pure Appl. Algebra 10 (1977), 271-294.
  
  \bibitem{Morrow}
M. Morrow, \emph{A case of the deformational Hodge conjecture via a pro Hochschild-Kostant-Rosenberg theorem}, C. R. Math. Acad. Sci. Paris 352 (3) (2014), 173-177. 
  
  \bibitem{Mu}
D. Mumford, \emph{Rational equivalence of 0-cycles on surfaces}, J. Math. Kyoto Univ. 9 (1968), 195-204.
  
  
  
\bibitem{PatelRavi1}
D. Patel and G. V. Ravindra, \emph{Weak Lefschetz for Chow Groups: Infinitesimal Lifting}, Homology, Homotopy, and Applications 16 (2) (2014), 65-84.

\bibitem{PW}
C. Pedrini and C. Weibel, \emph{Some surfaces of general type for which Bloch's conjecture holds}, Recent advances in Hodge theory, 308-329, London Math. Soc. Lecture Note Ser., 427, Cambridge Univ. Press, Cambridge, 2016.

   \bibitem{Quillen}
D. Quillen, \emph {Higher algebraic $K$-Theory I}, Lecture Notes in Mathematics, 341 (1972), 85-147.





 
 
 

\bibitem{Ro1}
A. A. Roitman, \emph{The torsion of the group of zero-cycles modulo rational equivalence}, Ann. of Math. 111 (1980),  553-569.

\bibitem{Ro2}
A. A. Roitman, \emph{Rational equivalence of zero-dimensional cycles}, (Russian) Mat. Sb.(N.S) 89 (131) (1972), 569-585.


\bibitem{Sernesi}
E. Sernesi, \emph{Deformations of algebraic schemes}, Grundlehren der Mathematischen Wissenschaften [Fundamental Principles of Mathematical Sciences], 334, Springer-Verlag, Berlin, 2006, xii+339 pp. ISBN: 978-3-540-30608-5; 3-540-30608-0.

\bibitem{Soule}
 C. Soul\'e, \emph{Op\'erations en K-th\'eorie alg\'ebrique}, Canad. J. Math. 37 (1985), 488-550.

\bibitem{Stien1}
J. Stienstra, \emph{On the formal completion of the Chow group $CH^{2}(X)$ for a smooth projective surface in characteristic 0}, Nederl. Akad. Wetensch. Indag. Math. 45 (3) (1983), 361-382. 

\bibitem{Stien2}      
J. Stienstra, \emph{Cartier-Dieudonn{\'e} theory for Chow groups,} J. Reine Angew. Math. 355 (1985), 1-66.

\bibitem{TT}
R. W. Thomason and T. Trobaugh, \emph{Higher algebraic K-theory of schemes and of derived
categories}, In The Grothendieck Festschrift, Volume III, volume 88 of Progress in Math.,
pages 247-436. Birkh$\ddot{\mathrm{a}}$user, Boston, Bas\'el, Berlin, 1990.

 \bibitem{VanderKallen}
W. van der Kallen, \emph { Le $K_2$ des nombres duaux}, C. R. Math. Acad. Sci. Paris 273 (1971), 1204-1207.
  
  
   \bibitem{Vistoli}
A. Vistoli, \emph{The deformation theory of local complete intersections}, available on arXiv:alg-geom/9703008.
  
  
  
  \bibitem{V1}
C. Voisin, \emph{Variations de structure de Hodge et z\'ero-cycles sur les surfaces g\'en\'erales}, Math. Ann. 299 (1) (1994), 77-103.
  
  \bibitem{V2}
C. Voisin, \emph{The generalized Hodge and Bloch conjectures are equivalent for general complete intersections}, Ann. Sci. {\'E}cole Norm. Sup 46 (3) (2013), 449-475.

\bibitem{V3}
C. Voisin, \emph{Bloch's conjecture for Catanese and Barlow surfaces}, J. Differential Geom. 97 (1) (2014), 149-175.


  
  
  \bibitem{W1}
 C. Weibel, \emph{Pic is a contracted functor}, Invent. Math. 103 (1991), 351--377.
 

\bibitem{W2}
C. Weibel, \emph{An introduction to homological algebra}, Cambridge Studies in Advanced Mathematics, 38, Cambridge University Press, Cambridge, 1994, xiv+450 pp. ISBN: 0-521-43500-5; 0-521-55987-1.

\bibitem{W3}
C. Weibel, \emph{Cyclic homology for schemes}, Proc. Amer. Math. Soc. 124 (6) (1996), 1655-1662.

 \bibitem{W4}
 C. Weibel, \emph{The Hodge filtration and cyclic homology}, K-Theory 12 (2) (1997), 145-164. 


\bibitem{Y1}
S. Yang, \emph{K-theory, local cohomology and tangent spaces to Hilbert schemes}, Annals of K-theory 3 (4) (2018), 709-722. 

\bibitem{Y2}
S. Yang, \emph{On extending Soul\'e's variant of Bloch-Quillen identification}, Asian J. of Math. 23 (1) (2019), 49-70. 


\bibitem{Y3}
S. Yang, \emph{Deformation of K-theoretic cycles}, Asian J. of Math. 24 (2) (2020), 303-330.

\bibitem{Y4}
S. Yang, \emph{Chern character and obstructions to deforming cycles}, J. of Algebra 601 (2022), 54-71.

\end{thebibliography}
\end{document}